\documentclass{ws-ijbc}
\usepackage{amsmath,bbm}
\usepackage{graphicx}
\usepackage{subfigure,captcont}
\usepackage{enumitem}

\newcommand{\fv}{\mathbf{f}}
\newcommand{\Fv}{\mathbf{F}}
\newcommand{\RR}{\mathbbm{R}}
\newcommand{\rv}{\mathbf{r}}
\newcommand{\sv}{\mathbf{s}}
\newcommand{\xv}{\mathbf{x}}

\DeclareMathOperator{\lcm}{lcm}

\begin{document}


\markboth{W.M. Seiler and M. Sei\ss}{No Chaos in Dixon's System}

\title{No Chaos in Dixon's System} 
\author{Werner M. Seiler}
\address{Institut f\"ur Mathematik, Universit\"at Kassel, 34109
  Kassel, Germany, seiler@mathematik.uni-kassel.de}
\author{Matthias Sei\ss}
\address{Institut f\"ur Mathematik, Universit\"at Kassel, 34109
  Kassel, Germany, mseiss@mathematik.uni-kassel.de}

\maketitle


\begin{abstract}
  The so-called Dixon system is often cited as an example of a
  two-dimensional (continuous) dynamical system that exhibits chaotic
  behaviour, if its two parameters take their value in a certain domain.
  We provide first a rigorous proof that there is no chaos in Dixon's
  system.  Then we perform a complete bifurcation analysis of the system
  showing that the parameter space can be decomposed into sixteen different
  regions in each of which the system exhibits qualitatively the same
  behaviour.  In particular, we prove that in some regions two elliptic
  sectors with infinitely many homoclinic orbits exist.
\end{abstract}


\section{Introduction}

\citet{cdk:dynmod} derived a model for the magnetic field of neutron stars
consisting of three autonomous ordinary differential equations.  One of the
equations can be easily decoupled so that the core of the model is given by
a two-dimensional dynamical system depending on two positive real
parameters $a$ and $b$ (see (\ref{eq:cdk}) below) which we will call
\emph{CDK system}.  In numerical simulations of this system, Cummings
\emph{et al.} observed for $0<a,b<1$ a complicated dynamics (which matches
some astronomical observations) and suspected a chaotic behaviour.  In a
subsequent article \cite{dck:ccd2d}, they reported -- based on additional
numerical experiments -- a positive Lyapunov exponent as a further
indication of chaos.

Since the Poincar\'e-Bendixson Theorem excludes chaos in two-dimensional
continuous dynamical systems, these findings raised some interest.
\citet{dck:ccd2d} already stated that the CDK system was not a
counterexample, as the assumptions of the theorem were violated.  In fact,
they were always very careful in their formulations: they used the word
chaos only with quotation marks and noted already in their first article
that classical conditions for chaos were not satisfied.  In a further
article, \citet{ddd:pdd} completely avoided the use of this term and coined
instead the phrase ``piecewise deterministic dynamics''.  Nevertheless,
under the name Dixon's system (which is a bit unfair towards the other
authors of these works) their differential equations are repeatedly cited
as an example of a chaotic system (see e.\,g.\ \cite{jcs:ec}), although
\citet{ade:2dchaos} already clearly indicated that the seemingly chaotic
behaviour is simply due to the existence of two elliptic sectors with
infinitely many homoclinic orbits. We will briefly indicate that also in a
logarithmic variant of the CDK system proposed by \citet{jcs:ec} as a
further ``chaotic'' system in the plane no chaos appears.

In this article, we will first exhibit some shortcomings of the analysis in
\cite{dck:ccd2d}.  Indeed, their argument against the applicability of the
Poincar\'e-Bendixson Theorem is inconsistent and furthermore relies on a
simplified version of the theorem.  Using a geometric approach based on
singularity theory, we will then provide a rigorous proof that also a
refined version cannot be directly applied to the CDK system.  We will show
that one can derive a polynomial vector field possessing the same integral
curves as the CDK system with the origin as a degenerate stationary point.
This allows us a rigorous qualitative analysis of the solution behaviour in
the neighbourhood of the origin.

\citet{cdk:dynmod} distinguished in their numerical studies four different
regions in the parameter space.  We will show that a rigorous
classification leads actually to sixteen different regions, if one also
takes the behaviour at infinity into account.  It turns out that although
the two parameters in the system are from a physical point of view of
essentially the same nature, they play mathematically a very different
role.  In particular, the emergence of the homoclinic orbits is completely
controlled by one of them alone.

\section{The CDK System and the Poincar\'e-Bendixson Theorem}
\label{sec:cdkpb}

In their work on the magnetic field of neutron stars, \citet{cdk:dynmod}
derived the following planar dynamical system
\begin{equation}\label{eq:cdk}
  \dot x = \frac{xy}{x^{2}+y^{2}}-a x\,,\qquad
  \dot y = \frac{y^{2}}{x^{2}+y^{2}}-b y+b-1
\end{equation}
depending on two positive parameters $a,b>0$ as decoupled part of a
three-dimensional model.  For the physical interpretation of the model, we
refer to their article.  The analysis of the stationary points of this
planar system and of their stability leads naturally to case distinctions
whether zero, one or two of the parameters have a value larger than $1$.
In numerical simulations, the case $0<a,b<1$ showed a complicated behaviour
that \citet{cdk:dynmod} suspected to be chaotic.

In the subsequent article, \citet{dck:ccd2d} provided a closer analysis of
this case giving further numerical evidence for a highly irregular dynamics
resembling chaos.  They also discussed whether their model represented a
counterexample to the Poincar\'e-Bendixson Theorem.  Obviously, the right
hand side of (\ref{eq:cdk}) consists of rational functions with a pole at
the origin.  \citet{dck:ccd2d} showed that its gradient does not possess a
well-defined limit for $(x,y)\rightarrow(0,0)$ and concluded that therefore
the system was not $\mathcal{C}^{1}$ as assumed in the Poincar\'e-Bendixson
Theorem.  This argument is erroneous in at least two respects.  Firstly,
the cited version of the Poincar\'e-Bendixson Theorem assumes a dynamical
system defined on the whole plane, whereas the right hand side of
(\ref{eq:cdk}) is obviously not defined at the origin.  Secondly, it makes
no sense to study the differentiability of a function at a point where it
is not defined.\footnote{The same erronous argument appears in
  \cite{ddd:pdd} where the Lipschitz continuity of a vector field is
  studied at a point where it is not defined.}  In fact, it is easy to
verify that the right hand side of (\ref{eq:cdk}) itself does not possess a
well defined limit for $(x,y)\rightarrow(0,0)$.

We will use here the following version of the Poincar\'e-Bendixson Theorem,
as it can e.\,g.\  be found in the text book by
\citet[Thm.~1, Sect.~3.7]{lp:deds}.

\begin{theorem}[Poincar\'e-Bendixson]\label{thm:pb}
  Let the autonomous dynamical system $\dot{\xv}=\fv(\xv)$ be defined on an
  open subset $\Omega\subseteq\RR^{2}$ with
  $\fv\in\mathcal{C}^{1}(\Omega,\RR^{2})$.  Assume that a forward
  trajectory lies completely in a compact subset
  $\mathcal{K}\subset\Omega$.  Then the $\omega$-limit set of the
  trajectory contains either stationary points of the system or it is a
  periodic orbit.
\end{theorem}

In the case of the CDK system (\ref{eq:cdk}), the obvious choice for the
open set $\Omega$ is the punctured plane $\Omega=\RR^{2}\setminus\{0\}$.
The decisive question for the applicability of Theorem~\ref{thm:pb} is now
whether we can find an enclosing compact set $\mathcal{K}\subset\Omega$ for
every bounded forward trajectory.  Only if this is the case, we can exclude
chaos by applying Theorem \ref{thm:pb} directly to (\ref{eq:cdk}).
However, we will show in Section \ref{sec:geoql} as a by-product of our
geometric analysis of (\ref{eq:cdk}) that this is not possible if the
values of the parameters $(a,b)$ lie in certain regions of the parameter
space.  More precisely, we will prove that for such parameter values there
exist infinitely many trajectories starting arbitrarily close to the origin
and approaching for $t\rightarrow\infty$ the origin again arbitrarily
close.  Obviously, any compact set $\mathcal{K}$ enclosing such a
trajectory must contain the origin and thus cannot be a subset of $\Omega$.
On the other hand, it is also clear that none of these trajectories exhibit
a chaotic behaviour, as each of them has an $\omega$-limit set consisting
simply of the origin.

However, even in this refined form this argument only shows that Theorem
\ref{thm:pb} is not directly applicable to the CDK system in the form
(\ref{eq:cdk}).  The key idea of our geometric analysis in Section
\ref{sec:geoql} will be to study instead of the rational system
(\ref{eq:cdk}) a polynomial vector field which is defined on the entire
plane $\RR^{2}$ and which is trajectory equivalent to (\ref{eq:cdk}) on
$\Omega$.  The application of Theorem \ref{thm:pb} to this system is
straightforward and shows immediately that the CDK system cannot exhibit
chaos.  The origin will be a degenerate stationary point for this
polynomial vector field and thus it becomes a standard problem in dynamical
systems theory to study the behaviour of solutions of (\ref{eq:cdk}) close
to the origin.

\section{A Geometric Approach to the CDK System}
\label{sec:geoql}

We follow here an approach pioneered by \citet{ves:int} who showed how the
analysis of an arbitrary system of differential equations can be reduced to
the study of (systems of) vector fields (see \cite{wms:invol} and
references therein for a modern representation).  In the case of ordinary
differential equations, this approach allows to transform implicit
equations into explicit ones.  In the context of the singularity theory of
ordinary differential equations -- see \cite{via:geoode} for a simple
introduction -- this approach is standard (however, without being
attributed to Vessiot).  In the literature, the arising theory is often
named after Cartan.  Following \citet{edf:vessiot}, we think that this is
not appropriate, as Cartan worked almost exclusively with differential
forms and not with vector fields.  Quasi-linear equations show in this
context a special behaviour with properties not present in arbitrary
non-linear equations.  Although this phenomenon has a very simple
explanation, it seems that it was exhibited for the first time only fairly
recently in \cite{wms:singbif}.  For a more detailed analysis see the
recent work \cite{wms:quasilin}.

In the geometric theory of differential equations (see \cite{wms:invol} and
references therein for an in-depth introduction into at least the regular
case), a differential equation is modelled as a fibred submanifold of an
appropriate jet bundle.  To be able to study also singularities, we take a
more general point of view and define a differential equation as a
subvariety of the jet bundle such that its image under the canonical
projection to the space of independent variables is dense.

More precisely, we take for the purpose of this article which studies
exclusively first-order ordinary differential equations with polynomial
nonlinearities as underlying base space $\RR$ and consider jets of local
functions $\RR\rightarrow\RR^{n}$.  Local coordinates of the first jet
bundle $\mathcal{J}_{1}(\RR,\RR^{n})\cong\RR\times\RR^{n}\times\RR^{n}$ are
then denoted by $(t,\xv,\dot{\xv})$.  We have two canonical projections
$\pi^{1}_{0}:\mathcal{J}_{1}(\RR,\RR^{n})\rightarrow\RR\times\RR^{n}$ and
$\pi^{1}:\mathcal{J}_{1}(\RR,\RR^{n})\rightarrow\RR$.  An \emph{algebraic
  differential equation} is a subvariety
$\mathcal{R}_{1}\subseteq\mathcal{J}_{1}(\RR,\RR^{n})$ such that
$\pi^{1}(\mathcal{R}_{1})$ lies dense in $\RR$.  Thus we consider
exclusively systems of the form $\Fv(t,\xv,\dot{\xv})=0$ where the
components of $\Fv$ are polynomials.

We distinguish two kinds of singularities of such a differential equation.
An \emph{algebraic singularity} is a point $\rho\in\mathcal{R}_{1}$ which
is not smooth, i.\,e.\ a singular point in the sense of algebraic geometry.
A \emph{geometric singularity} is a smooth point $\rho\in\mathcal{R}_{1}$
which is critical for the restricted projection
$\pi^{1}|_{\mathcal{R}_{1}}$, i.\,e.\ a singular point in the sense of
differential topology.  While much is known about the local solution
behaviour around geometric singularities (see e.\,g.\ \cite{via:geoode} for
an elementary introduction to the case of scalar first-order equations),
not much can be found in the literature about algebraic singularities.

The \emph{contact distribution} $\mathcal{C}^{(1)}$ of
$\mathcal{J}_{1}(\RR,\RR^{n})$ is generated by the $n+1$ vector fields
\begin{equation}\label{eq:cdq}
  C_{\mathrm{trans}}=\partial_{t}+
  \sum_{i=1}^{n}\dot{x}_{i}\partial_{x_{i}}\,,
  \qquad C_{i}=\partial_{\dot{x}_{i}}\,.
\end{equation}
The first field is transversal with respect to the projection $\pi^{1}$
(and essentially encodes the chain rule), whereas the remaining $n$ fields
$C_{i}$ span the vertical bundle for the projection $\pi^{1}_{0}$.  Given
an algebraic differential equation
$\mathcal{R}_{1}\subseteq\mathcal{J}_{1}(\RR,\RR^{n})$, we define the
\emph{Vessiot space} at a smooth point $\rho\in\mathcal{R}_{1}$ as
$\mathcal{V}_{\rho}[\mathcal{R}_{1}]=
T_{\rho}\mathcal{R}_{1}\cap\mathcal{C}^{(1)}_{\rho}$, i.\,e.\ as that part
of the contact distribution that is tangential to $\mathcal{R}_{1}$.
$\mathcal{V}_{\rho}[\mathcal{R}_{1}]$ is obviously a linear space whose
dimension can vary with the point $\rho$.  It is easy to see that on a
(Zariski) open subset of $\mathcal{R}_{1}$ the Vessiot spaces define a
smooth regular distribution.  By a certain abuse of language, we call the
whole family of Vessiot spaces the \emph{Vessiot distribution} of
$\mathcal{R}_{1}$.

We call an integral curve $\mathcal{S}\subseteq\mathcal{R}_{1}$ of the
Vessiot distribution $\mathcal{V}[\mathcal{R}_{1}]$ a \emph{generalised
  solution} of $\mathcal{R}_{1}$ and its projection
$\pi^{1}_{0}(\mathcal{S})$ a \emph{geometric solution}.  If the function
$\sv(t)$ is a solution of $\mathcal{R}_{1}$ in the classical sense, then
its graph is a geometric solution and the curve defined by all points
$\bigl(t,\sv(t),\dot{\sv}(t)\bigr)$ a generalised solution.  Note, however,
that not every geometric solution is necessarily the graph of a function.

Consider now a \emph{quasi-linear} system of the form
$A(t,\xv)\dot{\xv}=\rv(t,\xv)$ where we assume for simplicity that the
matrix $A$ is almost everywhere non-singular and that at those points where
it is singular the rank drops only by one.  Then it is shown in
\cite{wms:singbif} that outside of certain singularities the generator of
the Vessiot distribution is projectable to $\RR\times\RR^{n}$ and yields
there the vector field
$Y=\delta(t,\xv)\partial_{t}+
\bigl(C(t,\xv)\rv(t,\xv)\bigr)^{t}\partial_{\xv}$ where
$\delta(t,\xv)=\det{A(t,\xv)}$ and $C(t,\xv)$ is the adjunct of $A(t,\xv)$.
Strictly speaking, the projected vector field $Y$ is defined only on the
subset $\pi^{1}_{0}(\mathcal{R}_{1})\subseteq\RR\times\RR^{n}$.  However,
from the above explicit expression one can see that $Y$ can be analytically
extended to all points on $\RR\times\RR^{n}$ where $\delta$ and $C$ are
defined.  If we assume that $A$ and $\rv$ are polynomial in their
arguments, then obviously this means that we may consider $Y$ as a vector
field on the whole space $\RR\times\RR^{n}$.

We call an invariant curve of the vector field $Y$ a \emph{weak geometric
  solution}.  Indeed, any geometric solution is also a weak one; however,
parts of weak geometric solution may correspond to the graph of a function
which is continuous but not differentiable.  In particular, any integral
curve of the quasi-linear system is also an invariant curve of this vector
field.  Thus the analysis of the original implicit problem can be reduced
to the study of the vector field $Y$, i.\,e.\ an explicit autonomous
differential equation. For more details about this geometric approach we
refer again to \cite{wms:quasilin}.

For the application of the above sketched ideas to the CDK system
(\ref{eq:cdk}), we work in the jet bundle $\mathcal{J}_{1}(\RR,\RR^{2})$
with local coordinates $(t,x,y,\dot x,\dot y)$.  Its contact distribution
is locally generated by the following three vector fields
\begin{equation}\label{eq:cd}
  C_{\mathrm{trans}}=\partial_{t}+\dot x\partial_{x}+\dot y\partial_{y}\,,
  \qquad C_{x}=\partial_{\dot x}\,,\qquad C_{y}=\partial_{\dot y}\,.
\end{equation}

Instead of working with (\ref{eq:cdk}), we prefer to rewrite the system in
the quasi-linear form
\begin{equation}\label{eq:cdki}
  \begin{aligned}
    (x^{2}+y^{2})\dot x &= xy - a(x^{3}+xy^{2})\,,\\
    (x^{2}+y^{2})\dot y &= y^{2} - (b y-b+1)(x^{2}+y^{2})\,,
  \end{aligned}
\end{equation}
as these equations are also defined at points with $x=y=0$, since they are
polynomial in all variables.  Nevertheless, even in this implict form these
points are problematic.  As one can easily check with the Jacobian of
(\ref{eq:cdki}), the zero set of (\ref{eq:cdki}) is not a submanifold of
$\mathcal{J}_{1}(\RR,\RR^{2})$, but a three-dimensional variety for which
all these points represent a singularity in the sense of algebraic
geometry.

For the construction of the Vessiot distribution of the differential
equation $\mathcal{R}_{1}$ defined by (\ref{eq:cdki}), we take a general
vector field $X=\alpha C_{\mathrm{trans}}+\beta C_{x}+\gamma C_{y}$ in the
contact distribution $\mathcal{C}^{(1)}$ with yet undetermined coefficients
and check when it is tangential to $\mathcal{R}_{1}$.  This yields the
equations
\begin{equation}\label{eq:vess}
  \begin{aligned}
    (x^{2}+y^{2})\beta &= \bigl((y-3a x^{2}-a y^{2})\dot x+
                      (x-2a xy)\dot y\bigr)\alpha\,,\\
    (x^{2}+y^{2})\gamma &= -\bigl(2(b y-b+1)x\dot x+
                       b(3y+2)y\dot y\bigr)\alpha\,.
  \end{aligned}
\end{equation}

For the projected vector field $Y=(\pi^{1}_{0})_{*}(X)$ we only need a
value for $\alpha$ which obviously may be chosen as $x^{2}+y^{2}$.  Using
the equations (\ref{eq:cdki}), we thus obtain
\begin{equation}\label{eq:Y}
  Y=(x^{2}+y^{2})\partial_{t} +
    \bigl(xy - a(x^{3}+xy^{2})\bigr)\partial_{x} +
    \bigl(y^{2} - (b y-b+1)(x^{2}+y^{2})\bigr)\partial_{y}\,.
\end{equation}
Note that, strictly speaking, the Vessiot spaces are defined only at smooth
points, i.e.  in our example the vector field $X$ is not defined at all
points with $x=y=0$, as these are singularities in the sense of algebraic
geometry.  However, it is easy to see that the projected vector field $Y$
can be analytically extended to the whole of $\RR^3$, so that we can ignore
the fact that these points are algebraic singularities.

For analysing the phase portrait of (\ref{eq:cdki}), the dynamics in
$t$-direction is irrelevant, as it simply describes a reparametrisation of time, 
and we can concentrate on the planar dynamical
system defined by the second and third component of $Y$
\begin{equation}\label{eq:dynsys}
    \dot{x} = xy - a(x^{3}+xy^{2})\,,\quad
    \dot{y} = y^{2} - (b y-b+1)(x^{2}+y^{2})\,,
\end{equation}
where the dot denotes now the differentiation with respect to some parameter
$\tau$.  Alternatively, one may derive the same system via a time
reparametrisation in (\ref{eq:cdki}) and consider $\tau$ as the new time
variable.  We prefer our approach via singularity theory, as its
interpretation is more transparent in situations where the
reparametrisation does not define a bijective function.

Obviously, the polynomial dynamical system (\ref{eq:dynsys}) is defined on
the whole plane $\RR^{2}$ and outside of the origin its trajectories
coincide with the ones of the CDK system (\ref{eq:cdki}).  The
Poincar\'e-Bendixson Theorem even in its most elementary form trivially
applies to (\ref{eq:dynsys}) and thus we can rigorously exclude the
possibility of any form of chaotic dynamics in the CDK system.  It is easy
to see that the origin is a stationary point of (\ref{eq:dynsys}) and that
the Jacobian of the system vanishes at the origin.  Thus we are in the case
of a nilpotent stationary point requiring a careful analysis via blow-ups
to determine the local solution behaviour.

\section{Stationary Points of the CDK System}
\subsection{The Set of Stationary Points}
In this section we determine the set of stationary points of the dynamical 
system in dependence of the parameters $a, b \in \mathbb{R}_{>0}$. 
This set is given by the zero set of the polynomial system
\begin{subequations}
\label{eqn:condstatpoints}
 \begin{eqnarray}
  x (y-ay^2-ax^2) &=& 0
  \label{eqn:condstatpoint1}, \\
  -b y^3 - byx^2 +b y^2 + (b-1)x^2 &=& 0 
  \label{eqn:condstatpoint2}.
 \end{eqnarray}
\end{subequations}
In order to describe its solution set it is useful to distinguish between
the cases $x=0$ and $x\neq 0$.

For $x=0$, ($\ref{eqn:condstatpoint1}$) is
always satisfied and ($\ref{eqn:condstatpoint2}$) reduces to 
$-by^3 + b y^2 = b y^2 (1-y)=0$.
Solving this equation, we obtain the two zeros $(x,y)=(0,0)$ and $(x,y)=(0,1)$ 
which are independent from the parameter values. 
This means that the stationary points of the system with $x=0$ for all
$a$,$b \in \mathbb{R}_{>0}$ are $s_1 := (0,0)$ and $s_2 := (0,1)$. 

For $x\neq 0$, ($\ref{eqn:condstatpoint1}$) is satisfied if and only if  
\begin{equation}
\label{eqn:condstatpoint1xnonzero}
  x^2 = \frac{1}{a} y - y^2 
\end{equation}
and with this identity (\ref{eqn:condstatpoint2}) reduces to 
\begin{equation}
y ( \frac{a-b}{a} y + \frac{b-1}{a} ) =0 \label{eqn:cond1}.
\end{equation}
Equation~(\ref{eqn:cond1}) shows that the existence and the number of 
stationary points of the system depends on the values of the parameters. 
We obtain the following cases:
\begin{description}[nosep]
\item[(a)] If the parameter satisfy $b=1$, $ a \neq b$ or $b\neq 1$, $a=b$,
  then there are no stationary points with $x\neq 0$.  Indeed, for these
  values one of the two coefficients of (\ref{eqn:cond1}) vanishes and so
  (\ref{eqn:cond1}) has exactly one solution, namely $y=0$.  Reducing
  (\ref{eqn:condstatpoint1}) with this solution, we see that the only
  possible solution for it is $x=0$.
\item[(b)] If the parameters are $a=1$, $b=1$, then there are infinitely
  many stationary points with $x \neq 0$. More precisely, in this case both
  coefficients of (\ref{eqn:cond1}) vanish. Thus $y$ is arbitrary and we
  deduce from (\ref{eqn:condstatpoint1}) that the stationary points are the
  points of the circle $x^2+(y-\tfrac{1}{2})^2=\tfrac{1}{4}$ with
  $x \neq 0$. Note that by dropping the condition $x \neq 0$, we may
  include the stationary points $s_1$ and $s_2$ into the circle.
\item[(c)] If the parameters are $b \neq 1$, $a \neq b$, then
  (\ref{eqn:cond1}) has exactly two solutions, as both coefficients do not
  vanish.  It is easy to check that for these parameter values the
  solutions of (\ref{eqn:cond1}) are $y_1 :=0$ and
\begin{equation}
\label{eqn:statpointssol}
      y_2:=-\frac{b-1}{a-b}. 
\end{equation}
On the other hand, (\ref{eqn:condstatpoint1}) can only have a solution with
$x \neq 0$ if $y$ is also non-zero. So the only possibility for the
$y$-component of a solution of (\ref{eqn:condstatpoints}) is
$y_2$. Substituting $y_2$ for $y$ in (\ref{eqn:condstatpoint1xnonzero}) and
solving for $x^2$ yields
\begin{equation}
\label{eqn:cond1y2}
    x^2 = - \frac{(b-1)^2}{(a-b)^2} - \frac{b-1}{a(a-b)} .
\end{equation}
and this equation has two real solutions if and only if the negative of the
first term of its right hand side is smaller than the second term. Since
the real numbers $(b-1)^2$, $(a-b)^2$ and $a$ are positive, we conclude
that $(b-1)$ and $(a-b)$ must have opposite sign. So if we multiply the
condition for real roots with $(a-b)/(b-1)$ and then taking the reciprocal
we get $(a-b)/(b-1) < -a$. There are two possibilities for the left hand
side of this inequation to be negative, namely $b >\max{\{1,a\}}$ and
$b<\min{\{1,a\}}$. In the first case the inequation rewrites as
$a-b < -ab+a$ and so $1 > a$. In the second case we obtain $1<a$ from
$a-b > -ab +a$.
\end{description}
    
It is easily checked that the parameter values specified in \textbf{(a)},
\textbf{(b)} and \textbf{(c)} exhaust the positive quadrant of
$\mathbb{R}^2$, that is, we determined the stationary points for all
possible parameter values $a$, $b \in \mathbb{R}_{>0}$.

\begin{lemma}\label{lem:stationarypoints}
  Depending on $a,b>0$, the dynamical system (\ref{eq:dynsys}) has the
  following finite stationary points:
  \begin{enumerate}[itemindent=20pt,itemindent=20pt,nosep]
  \item For $a=b=1$, there are infinitely many finite stationary points,
    namely the points of the circle \label{lem:stationarypointsCase1}
    \begin{equation*}
      x^2+(y-\tfrac{1}{2})^2=\tfrac{1}{4}\,.
    \end{equation*}
  \item \label{lem:stationarypointsCase2} For $b>1>a$ and $b<1<a$, there
    are four finite stationary points, namely $s_1=(0,0)$, $s_2 =(0,1)$,
    $s_3=(x_1,y_2)$ and $s_4=(x_2,y_2)$ where $y_2$ is given by
    (\ref{eqn:statpointssol}) and
    \begin{equation}\label{eq:x12}
         x_{1/2}=\pm \sqrt{\frac{y_2}{a} - y_2^2}\,.
    \end{equation}
  \item \label{lem:stationarypointsCase3} In all other cases, the system
    has two finite stationary points, namely $s_1=(0,0)$ and $s_2 =(0,1)$.
\end{enumerate}
\end{lemma}

\subsection{The Dynamics at the Stationary Point $s_2$}

In this section we study the dynamics of the system at the stationary point
$s_2$ using standard techniques and results from dynamical system theory
(for a nice presentation see for example \cite{dla:qualplan}).  The
following analysis of the dynamics at $s_2$ will only be valid for those
parameter values $a$, $b \in \mathbb{R}_{>0}$ for which $s_2$ is an
isolated singularity, that is, we exclude the case $a=b=1$ from our
considerations.

In order to apply standard results, we move the stationary point $(0,1)$ to
the origin using the transformation $(x,y)=(\bar{x},\bar{y}+1)$. After this
shift the new system is
\begin{subequations}
\label{eqn:completeSystemPoint_s2}
\begin{eqnarray}
 \dot{\bar{x}} &=& \bar{x}(\bar{y}(1-2a-a \bar{y})-a \bar{x}^2+1-a ), 
 \label{eqn:system_s2_eq1} \\
 \dot{\bar{y}} &=& -b \bar{y}(\bar{y}^2 +2\bar{y}+\bar{x}^2+1)-\bar{x}^2. 
 \label{eqn:system_s2_eq2} 
\end{eqnarray}
\end{subequations}
The dynamical behaviour (\ref{eqn:completeSystemPoint_s2}) at the
stationary point $(0,0)$ is determined by the eigenvalues of its linear
part.  The Jacobian matrix of the system at the origin computes as
$\mathrm{diag}(-a+1,-b)$.  The values of the eigenvalues clearly depend on
the parameter values for which $s_2$ is an isolated stationary point.  We
distinguish the following cases:
\begin{description}[nosep]
\item[(a)] For $a < 1$ the eigenvalues are non-zero and so the stationary
  point is hyperbolic. Moreover, they have clearly opposite sign and so by
  \cite[ Theorem 2.15 (i)]{dla:qualplan} the stationary point is a saddle
  point where the points on the invariant analytic curve tangent to the
  $\bar{x}$-axis are repelled and the points on the invariant analytic
  curve tangent to the $\bar{y}$-axis are attracted to the origin.
\item[(b)] For $a = 1$ the first eigenvalue is zero and the second one is
  non-zero and negative, that is, the origin is a semi-hyperbolic
  singularity in this case. Applying \cite[Theorem 2.19]{dla:qualplan}, one
  readily checks that it is an attracting node for $b>1$ and a saddle for
  $b<0$ with the $\bar{x}$-direction repelling and the $\bar{y}$-direction
  attractive.
\item[(c)] For $a > 1$, the eigenvalues are non-zero and so the stationary
  point is hyperbolic.  Since the eigenvalues are negative we deduce with
  \cite[Theorem 2.15 (ii)]{dla:qualplan} that the origin in this case is an
  attracting node.
\end{description}
 
\begin{lemma}\label{lem:s2}
  The dynamics at $s_2$ in the situation of
  Lemma~\ref{lem:stationarypoints} (\ref{lem:stationarypointsCase2}) and
  (\ref{lem:stationarypointsCase3}) is as follows:
  \begin{enumerate}[itemindent=20pt,nosep]
  \item For $a<1$ the stationary point $s_2$ is a saddle point where the
    unstable manifold is tangent to the $x$-axis and the stable manifold to
    the $y$-axis.
  \item For $a=1$ and $b>1$ the stationary point $s_2$ is an attracting
    node and for $a=1$ and $b<1$ it is a saddle point where the unstable
    manifold is tangent to the $x$-direction and the stable manifold to the
    $y$-axis.
  \item For $a>1$ the stationary point $s_2$ is an attracting node.
  \end{enumerate}
\end{lemma}

\subsection{The Dynamics at the Stationary Point $s_1$}

In this section we determine the dynamics near the stationary point $s_1$.
As in the case of $s_2$, we consider only those parameter values $a$,
$b \in \mathbb{R}_{>0}$ for which $s_1$ is an isolated
singularity. Therefore throughout this section we assume that the
parameters have values as described in Lemma~\ref{lem:stationarypoints}
(\ref{lem:stationarypointsCase2}) or (\ref{lem:stationarypointsCase3}).

In case of a hyperbolic stationary point, the dynamical behaviour of a
system can be read off from the eigenvalues of its linear
part. Unfortunately, the stationary point $s_1$ of our system is
non-elementary, as the Jacobian at this point is identically zero. Indeed,
it is easily checked that the Jacobian computes as
\[
    J(x,y)=
    \begin{pmatrix}
    -a y^2+y-3 a x^2 &  x(1-2 a y) \\
    2 b x(1- y) -2x & by(2 -3  y)- bx^2  
    \end{pmatrix} 
\]
and thus $J(0,0)$ is the zero matrix.  

A basic tool for the determination of the dynamics at a non-hyperbolic
stationary point is blowing-up. A good overview can be found in
\cite[Section 3]{dla:qualplan} whose notation we will use throughout. For
the desingularisation of $s_1$ we will use quasi-homogeneous directional
blow-ups. The determination of these blow-ups will depend on the parameters
$a$, $ b \in \mathbb{R}_{>0}$.  It turns out that is is useful to
distinguish between the cases $b \neq 1$ and $b=1$.

We start with $b \neq 1$. In this case no term in the dynamical system
vanishes.  From the Newton polygon we obtain the exponents $(1,1)$ so that
we are actually performing a homogeneous blow-up.  For the blow-up in
positive $x$-direction, we use the transformation
$(x,y)= (\bar{x},\bar{x} \bar{y})$.  The blown-up system is given by
\begin{subequations}
\begin{eqnarray}
\dot{\bar{x}} \ &=& \ \bar{x} ( \bar{y}-a\bar{x}(1+\bar{y}^2)), 
\label{eqn:xdirectionblow3} \\
\dot{\bar{y}} \ &=& \ (b-1)(\bar{y}^2 +1) + \bar{x}\bar{y}^2
                    (a-b)(\bar{y}+1).
                    \label{eqn:xdirectionblow4} 
\end{eqnarray}
\end{subequations}
We determine the stationary points of the system on the line $\bar{x}=0$.
The conditions are given by setting the right hand side of
(\ref{eqn:xdirectionblow3}) and (\ref{eqn:xdirectionblow4}) to zero. For
$\bar{x}=0$ the first condition is for all values of $\bar{y}$ fulfilled
and the second condition reduces on $\bar{x}=0$ to
$(b-1) (\bar{y}^2 +1)=0$.  Since $b\neq 1$, the last equation can only be
zero if and only if $(\bar{y}^2 +1)=0$ and so there are no real stationary
points of the system on $\bar{x}=0 $. The blow-up in negative $x$-direction
provides no new information, because the weight $\alpha$ is odd.

For the blow-up in positive $y$-direction, we consider the change of
variables $(x,y)= ( \bar{x} \bar{y},\bar{y})$ and obtain
\begin{eqnarray*}
  \dot{\bar{x}} \ &=&
  \ (\bar{x}^2+1) \big( \bar{y} \bar{x} (b-a) + \bar{x} (1-b)\big), \\
  \dot{\bar{y}} \ &=&
  \ \bar{y} \big( b+(b-1)\bar{x}^2 - b \bar{y}(\bar{x}^2+1)\big). 
\end{eqnarray*}
On the line $\bar{y}=0$, the conditions for a stationary point reduce to
the single equation $\bar{x} (\bar{x}^2+ 1)(1-b)=0$ with $\bar{x}=0$ as
only real solution.  Thus $(0,0)$ is the only stationary point of the
system with $\bar{y}=0$. We determine the dynamics near it by studying the
linear part of the system. It is easily checked that the Jacobian of the
system reduces at $(\bar{x},\bar{y})=(0,0)$ to
$J(0,0)=\mathrm{diag}(1-b,b)$.  Since by assumption $b\neq 1 $, the point
$(0,0)$ is an elementary hyperbolic stationary point. Moreover, if $b<1$,
then the eigenvalues are positive, and if $b>1$, then one eigenvalue is
negative and the other one is positive.  By \cite[Theorem
2.15]{dla:qualplan}, in the first case the stationary point is a repelling
node and that in the second case it is a saddle point.  In the latter case,
the points on the invariant analytic curve tangent to the $\bar{x}$-axis
are attracted towards the origin and on the curve tangent to the
$\bar{y}$-axis they are repelled from the origin.

Since $\beta$ is odd, the blow-up in negative $y$-direction yields similar
results.  As in the latter case, the only stationary point of the new
system is $(0,0)$.  One checks that the Jacobian matrix of the system
reduces at this point to $J(0,0)= \mathrm{diag}(-1+b,-b)$.  The assumptions
on $b$ imply that the eigenvalues are non-zero and so the stationary point
is again elementary hyperbolic. For $b<1$ the eigenvalues are obviously
negative and for $b>1$ the eigenvalues are distinct and of opposite
sign. Again by \cite[Theorem 2.15]{dla:qualplan} the stationary point is an
attracting node in the first case and in the second case it is a saddle
point where compared to the blow-up in positive $y$-direction the
orientation of trajectories is now inverted.

\begin{figure}[ht]
    \centering
    \begin{minipage}{0.5\textwidth}
      \includegraphics[width=0.45\textwidth]{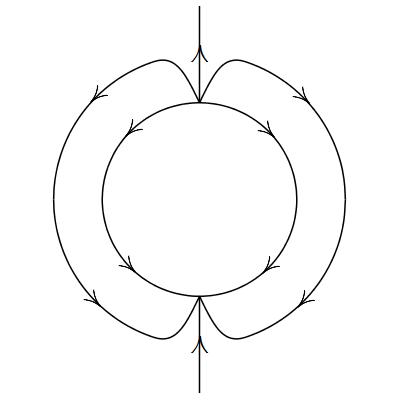}%
      \includegraphics[width=0.45\textwidth]{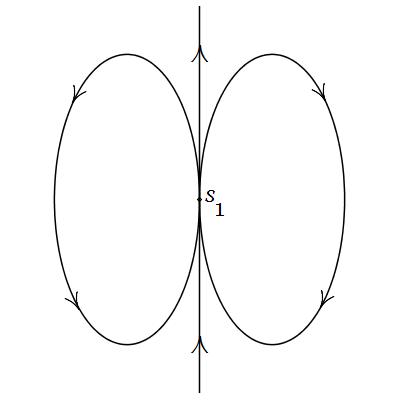}
    \caption{Blow-up and blow-down for $b<1$}
    \label{BlowUpforb<1}
    \end{minipage}%
    \begin{minipage}{0.5\textwidth}
      \includegraphics[width=0.45\textwidth]{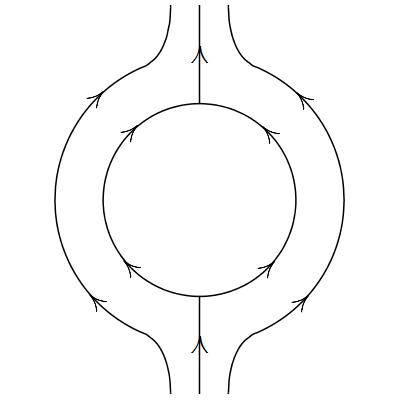}%
      \includegraphics[width=0.45\textwidth]{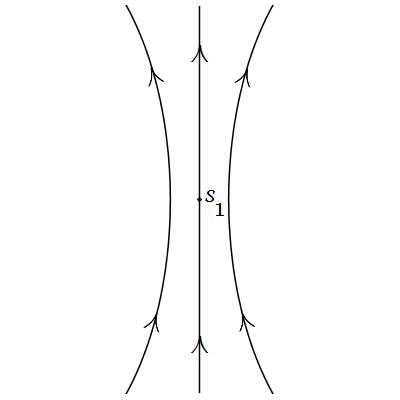}
    \caption{Blow-up and blow-down for $b>1$}
    \label{BlowUpforb>1}
    \end{minipage}
\end{figure}

In Figs.~\ref{BlowUpforb<1} and \ref{BlowUpforb>1} the discussed blow-ups
and their blow-downs are represented graphically. They show the local phase
portraits near $s_1$ in dependence of the parameter values $b<1$ and $b>1$.
\noindent

We consider now the case $b=1$. In this situation the term $(b-1)x^2$ in
the second equation of \eqref{eq:dynsys} vanishes. This has the effect
that, compared to the case $b\neq 0$, the Newton polygon yields different
weights $(\alpha,\beta)$ for the quasi-homogeneous directional blow-up,
namely $(1,2)$.

First we perform the blow-up in positive $x$-direction. With the above
weights, the change of variables for the positive $x$-direction is
$(x,y)=( \bar{x},\bar{x}^2 \bar{y})$. The blown-up vector field is
\begin{eqnarray*}
  \dot{\bar{x}} \ &=& \ \bar{x} ( \bar{y}-a-a\bar{x}^2 \bar{y}^2 ),  \\ 
  \dot{\bar{y}} \ &=& \  \bar{y} ( (2a-1) (\bar{x}^2 \bar{y}^2 +1) -\bar{y} ).
\end{eqnarray*}
We determine its stationary points on $\bar{x}=0$. Clearly the right hand
side of the first equation always vanishes on $\bar{x}=0$.  Reducing the
right hand side of the second equation with $\bar{x}=0$ leaves us with the
condition $\bar{y} ( (2a-1) - \bar{y} ) =0$ for a stationary point.  This
last equation is satisfied if and only if $\bar{y}=0$ or $\bar{y}=2a-1$.
Thus on $\bar{x}=0$ the stationary points of the system are $(0,0)$ and
$(0,2a-1)$.  The dynamics in a neighborhood of these points can be read off
from the eigenvalues of their linear parts. The Jacobian of the system
reduces at these points to $J(0,0) = \mathrm{diag}(-a, 2a-1)$ and
$J(0,2a-1) = \mathrm{diag} (a-1,1-2a)$.

Obviously, the eigenvalues of these Jacobians depend on the parameter
$a>0$.  The values of $a$ for which an eigenvalue reduces to zero are
clearly $a=1/2$ and $a=1$. Recall that we only consider in this section
parameter values such that $s_1$ is an isolated stationary point, that is
we ignore the case $a=1$.  Overall we distinguish between the following
cases:
\begin{description}[nosep]
\item[(a)] If $a< 1/2$, then the first Jacobian matrix has two negative
  eigenvalues and the second one has a negative and a positive
  eigenvalue. We conclude with \cite[Theorem 2.15]{dla:qualplan} that the
  stationary point $(0,0)$ is an hyperbolic attracting node and that the
  stationary point $(0,2a-1)$ is a saddle.  Note that in the latter case
  the point lies in the half plane $ y \leq 0$.  The phase portrait of this
  blow up is represented graphically in the right half sphere of
  Figure~\ref{BlowUpforb=1a<05} with corresponding stationary points $p_1$
  and $p_2$.
\item[(b)] If $a = 1/2$, then the two stationary points and their Jacobian
  matrices coincide. One eigenvalue is negative and the other one is
  zero. So the stationary point $(0,0)$ is semi-hyperbolic. In order to
  apply \cite[Theorem 2.19]{dla:qualplan} we consider the negative vector
  field. Then the nonzero eigenvalue is positive and we need to find a
  solution $x=f(y)$ in a neighbourhood of $(0,0)$ of the equation
  $\frac{1}{2}x^3y^2-xy+\frac{1}{2}x=0$ where the notation is as in the
  theorem.  A solution is clearly $x=f(y)=0$ and substituting it into
  $g(x)=y^2$, the theorem implies that the stationary point is a
  saddle-node. The graphic representation of this blow up is given in the
  right half plane of Figure~\ref{BlowUpforb=1a=05} where the stationary
  point is denoted by $p_1$.
\begin{figure}[ht]
    \centering
    \begin{minipage}{0.5\textwidth}
      \includegraphics[width=0.45\textwidth]{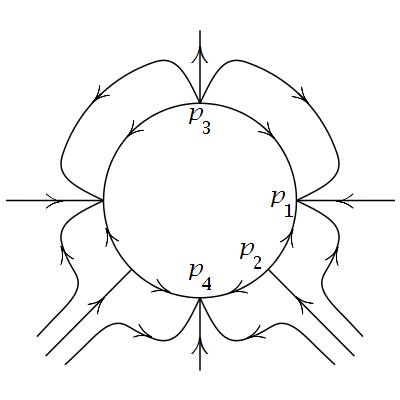}%
      \includegraphics[width=0.45\textwidth]{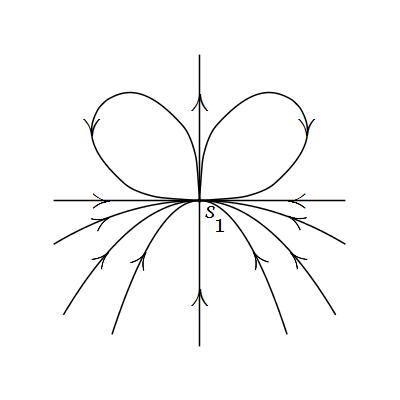}
    \caption{Blow-up and blow-down for $b=1$ and $a<0.5$}
    \label{BlowUpforb=1a<05}
    \end{minipage}%
    \begin{minipage}{0.5\textwidth}
      \includegraphics[width=0.45\textwidth]{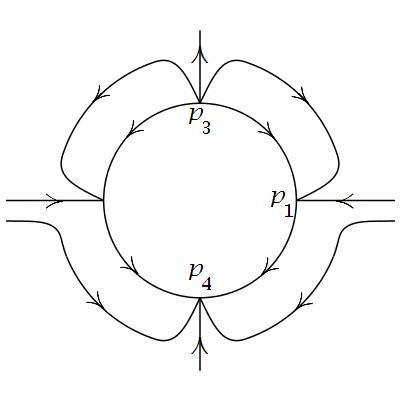}%
      \includegraphics[width=0.45\textwidth]{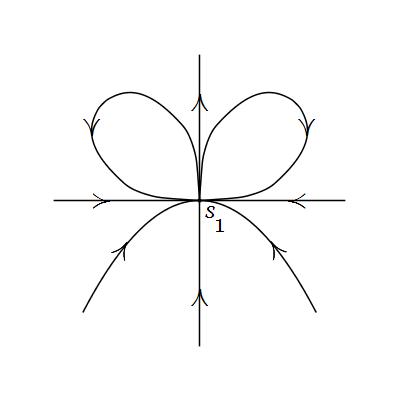}
    \caption{Blow-up and blow-down for $b=1$ and $a=0.5$}
    \label{BlowUpforb=1a=05}
    \end{minipage}
\end{figure}

\begin{figure}[ht]
    \centering
    \begin{minipage}{0.5\textwidth}
      \includegraphics[width=0.45\textwidth]{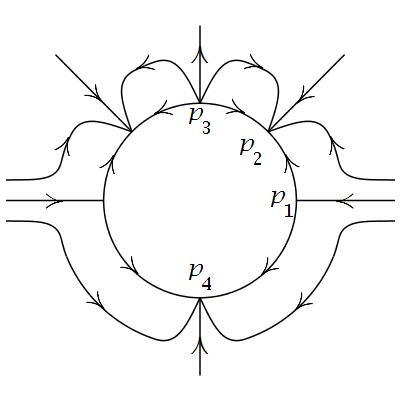}%
      \includegraphics[width=0.45\textwidth]{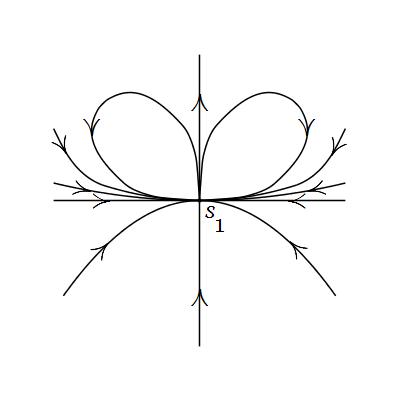}
    \caption{Blow-up and blow-down for $b=1$ and $0.5<a<1$}
    \label{BlowUpforb=105<a<1}
    \end{minipage}%
    \begin{minipage}{0.5\textwidth}
      \includegraphics[width=0.45\textwidth]{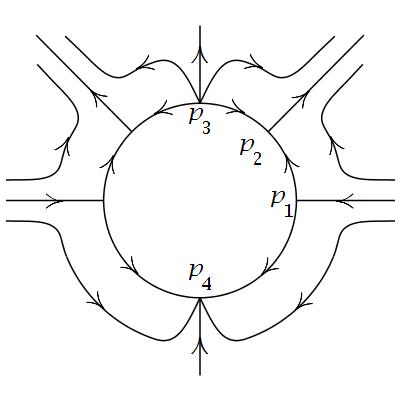}%
      \includegraphics[width=0.45\textwidth]{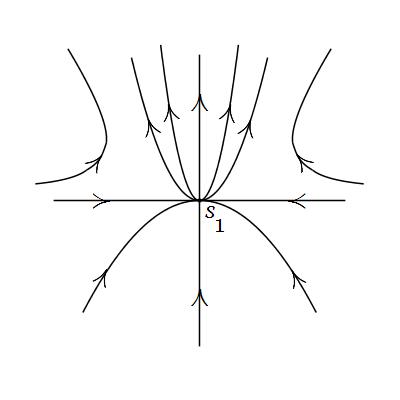}
    \caption{Blow-up and blow-down for $b=1$ and $a>1$}
    \label{BlowUpforb=1a>1}
    \end{minipage}
\end{figure}

\item[(c)] If $1 > a> 1/2$, then the two eigenvalues of $J(0,0)$ have
  opposite sign and the eigenvalues of $J(0,2a-1)$ are both negative, that
  is, both stationary points are hyperbolic.  We conclude with
  \cite[Theorem 2.15]{dla:qualplan} that $(0,0)$ is a saddle point and that
  $(0,2a-1)$ is an attracting node. Note that in this case the point
  $(0,2a-1)$ lies now in the half plane $ y \geq 0 $.  This results of this
  discussion are represented graphically in the right half plane of
  Figure~\ref{BlowUpforb=105<a<1} where $p_1$ and $p_2$ stand for $(0,0)$
  and $(0,2a-1)$.
\item[(d)] If $a> 1$, then the two eigenvalues of both matrices have
  opposite sign where $J(0,0)$ has a negative and a positive eigenvalue and
  $J(0,2a-1)$ has a positive and negative eigenvalue. So by \cite[Theorem
  2.15]{dla:qualplan} both points are saddle points. In the first case the
  direction of the $x$-axis is repelling and the direction of the $y$-axis
  is attracting. In the second case we observe exactly the opposite
  behaviour. Note that the point $(0,2a-1)$ lies again in the half plane
  $y \geq 0$. The phase portrait of this blow-up is sketched in the right
  half plane of Figure~\ref{BlowUpforb=1a>1}, where the stationary points
  are represented by $p_1$ and $p_2$ correspondingly.
\end{description}

The blow-up in negative $x$-direction yields no new information, since
$\alpha$ is odd. According to the same case distinctions, we obtain the
same stationary points with the same dynamical behaviour.  It is sketched
in the left half plane of
Figs.~\ref{BlowUpforb=1a<05}--\ref{BlowUpforb=1a>1}.
  
For the blow-up in positive $y$-direction with weights $(1,2)$, we use the 
transformation $(x,y)=( \bar{x} \bar{y}, \bar{y}^2)$ leading to the system 
\begin{eqnarray*}
  \dot{\bar{x}}&=& - \bar{x} \left( (2 a-1)(\bar{y}^2 +\bar{x}^2 ) -1 \right), \\
 \dot{\bar{y}}&=&  - \bar{y} (\bar{y}^2 + \bar{x}^2-1). 
\end{eqnarray*}
We determine the stationary points of the system on $\bar{y}=0$.  Obviously
for $\bar{y}=0$ the right hand side of the second equation vanishes and so
the condition for a stationary point on $\bar{y}=0$ is given by
$\bar{x} ((2a-1)\bar{x}^2-1)=0 $.  We conclude that the stationary points
on $\bar{y}=0$ are $(0,0)$ and $(\pm 1/\sqrt{2a-1},0)$. The latter two
points are real only for $a>\frac{1}{2}$ and we already analysed these two
points in the blow-up in positive $x$-direction in the respective case. So
it is left to determine the qualitative behaviour near $(0,0)$. It is
easily checked that the Jacobian reduces at $\bar{x}=0$ and $\bar{y}= 0$ to
the unit matrix.  This means that independent of the parameter values we
have two positive eigenvalues and so by \cite[Theorem 2.15]{dla:qualplan}
the stationary point is a repelling node.  The results of this discussion
are represented graphically in the upper half plane of
Figs.~\ref{BlowUpforb=1a<05}--\ref{BlowUpforb=1a>1} where the stationary
point is denoted by $p_3$.

For the blow-up in negative $y$-direction with the same weights we use
$(x,y)=( \bar{x} \bar{y}, -\bar{y}^2)$ and obtain
\begin{eqnarray*}
\dot{\bar{x}}&=& - \bar{x} \left( (2a-1) (\bar{y}^2 + \bar{x}^2 ) + 1 \right), \\
\dot{\bar{y}}&=& - \bar{y} (\bar{y}^2+ \bar{x}^2+1). 
\end{eqnarray*}
Similar as above, we conclude that the condition for a stationary point of
the system on $\bar{y}=0$ is given by the equation
$\bar{x} ((2a-1) \bar{x}^2 + 1)=0$.  Clearly $(0,0)$ is a stationary point
independent of the parameter value $a>0$.  For $a< \frac{1}{2}$ the
coefficient $(2a-1)$ is negative and so there are two additional real
stationary points, namely $(\pm 1/\sqrt{ -(2a-1)}),0)$.  The dynamics of
these two points have already been discussed in the blow-ups in positive
and negative $x$-direction. So we only need to analyse the dynamics near
$(0,0)$.  It is easily checked that the Jacobian reduces at
$(\bar{x},\bar{y})=(0,0)$ to the negative of the unit matrix and so the
stationary point $(0,0)$ is an attracting node by \cite[Theorem
2.15]{dla:qualplan}. In
Figs.~\ref{BlowUpforb=1a<05}--\ref{BlowUpforb=1a>1}, this point is called
$p_4$.

\begin{lemma}\label{lem:stationaryPoint_s1}
  Assume we are in case (\ref{lem:stationarypointsCase2}) or
  (\ref{lem:stationarypointsCase3}) of Lemma~\ref{lem:stationarypoints}.
  Then at the isolated stationary point $s_1$ we have the following
  dynamical behaviour:
  \begin{enumerate}[itemindent=20pt,nosep]
  \item For $b<1$ there are two elliptic sectors as in
    Fig.~\ref{BlowUpforb<1}.
  \item For $b>1$ there are two parabolic sectors as in
    Fig.~\ref{BlowUpforb>1}.
  \item For $b=1$ there are four subcases:
    \begin{enumerate}[leftmargin=20pt,itemindent=15pt,nosep]
    \item For $a < 1/2$ there are two elliptic sectors in the upper half
      plane and four hyperbolic sectors in the lower half plane as in
      Fig.~\ref{BlowUpforb=1a<05}.
    \item For $a=1/2$ there are two elliptic sectors in the upper half
      plane and two hyperbolic sectors in the lower half plane as in
      Fig.~\ref{BlowUpforb=1a=05}.
    \item For $1>a>1/2$ there are in the upper half plane two elliptic and
      two hyperbolic sectors and two additional hyperbolic sectors in the
      lower half plane as in Fig.~\ref{BlowUpforb=105<a<1}.
    \item For $a>1$ there are in the upper half plane two parabolic and two
      hyperbolic sectors and two additional hyperbolic sectors in the lower
      half plane as in Fig.~\ref{BlowUpforb=1a>1}.
    \end{enumerate}
  \end{enumerate}
\end{lemma}

\subsection{The Dynamics at the Stationary Points $s_3$ and $s_4$}

In this section we determine the dynamics near the stationary points $s_3$
and $s_4$.  Recall that, according to Lemma~\ref{lem:stationarypoints},
these points are stationary if and only if the parameters $a$, $b>0$
satisfy either $b>1>a$ or $a>1>b$.  Therefore we assume throughout this
section that $a$ and $b$ are chosen such that one of these relations is
fulfilled.

In order to apply the standard results of \cite{dla:qualplan}, we need to
move each of the stationary points $s_3$ and $s_4$ with a correspond change
of variables to the origin.  By Lemma~\ref{lem:stationarypoints}, the
stationary points are $s_3=(x_1,y_2)$ and $s_4=(x_2,y_2)$. Thus the
transformation $(x,y)=(x_{k} + \bar{x},y_2 + \bar{y})$ moves for $k=1$ the
point $s_3$ and for $k=2$ the point $s_4$ to the origin.  Applying this
transformation to the original system, we obtain the new system
\begin{eqnarray*}
 \dot{\bar{x}} &=& - ( (\bar{x}^2 + 2 x_{k} \bar{x} + \bar{y}^2 + 2 y_2 \bar{y} 
 + x_{k}^2 + y_2^2 )a - \bar{y} - y_2 )( \bar{x}+x_{k}),  \\
 \dot{\bar{y}} &=& -(\bar{y}+y_2-1)(\bar{x}^2 + 2 x_{k} \bar{x} + \bar{y}^2 + 
 2 y_2 \bar{y} +x_{k}^2 +y_2^2)b -( \bar{x} + x_{k})^2 ,
\end{eqnarray*}
which depends obviously on the values for $x_{k}$ and $y_2$. In order to
determine its dynamical behaviour near the origin, we compute the
eigenvalues of its linear part.  If we reduce the entries of the Jacobian
matrix with $y_2=-(b-1)/(a-b)$ and $(\bar{x},\bar{y})=(0,0)$, we get the
matrix
\[
J(0,0)=
\begin{pmatrix}
\frac{2b(b-1)(a-1)}{(a-b)^2} & \frac{ ( (2b-1)a-b )  x_{k} }{a-b}  \\
\frac{2a(b-1)x_{k}}{a-b}     & -\frac{b(b-1)(2a^2-3a+b)}{a(a-b)^2}
\end{pmatrix}
\]
whose entries depend on the value for $x_{k}$. The eigenvalues of $J(0,0)$
compute as
\begin{equation*}
\lambda_{1/2} =\frac{ b (1-b) \pm (b-1)\sqrt{b ( b+ 8a(a-1))}}{2a (b-a)}\,. 
\end{equation*}
We see that they do not depend on $x_{k}$ and so the dynamics near the
stationary points $s_3$ and $s_4$ of the original system are qualitatively
the same. For its determination it is useful to make the following case
distinctions for the parameters $a$ and $b$:
\begin{description}[nosep]
\item[(a)] Suppose that $b>1>a$. It follows that the expression $8a(a-1)$
  is a negative real number. If we assume in addition that
  $| 8a(a-1) | > b$, then the imaginary parts of the eigenvalues
  $\lambda_{1/2}$ are nonzero and $\lambda_{1/2}$ are complex conjugated to
  each other. Further from $b>1$ and $b>a$ it follows that their common
  real part is negative. By \cite[Theorem 2.15 (iii)]{dla:qualplan} the
  origin is in this situation a strong attracting focus.  Now suppose that
  $| 8a(a-1) | \leq b$. Then both eigenvalues are real and in order to
  determine the dynamics we need to decide whether they are positive or
  negative. It follows from the assumption that $2a(b-a)$ is positive and
  so $b(1-b)/(2a(b-a))$ is negative. Because $b-1$ and
  $\sqrt{b ( b+ 8a(a-1))}$ are both positive, the second summand in the
  expression for the eigenvalues is also positive. We conclude that the
  eigenvalue $\lambda_2$ is negative.  The eigenvalue $\lambda_1$ is also
  negative. Indeed, the inequation $8a(a-1)<0$ implies that
  $\sqrt{b ( b+ 8a(a-1))} < b$ and so we have
  $b (1-b) + (b-1)\sqrt{b ( b+ 8a(a-1))} <0$.  Because both eigenvalues are
  real and negative we obtain from \cite[Theorem 2.15 (ii)]{dla:qualplan}
  that the stationary point is an attracting node.
\item[(b)] Assume that $a>1>b$. It is easily seen that then the eigenvalues
  are real numbers and we need to determine if they are positive are
  negative.  The real numbers $2a(b-a)$ and $b-1$ are clearly negative and
  the product $b(1-b)$ is positive. We conclude that $b(1-b)/(2a(b-a))<0$
  and $(b-1)/(2a(b-a))>0$.  It can be easily checked now that $\lambda_2$
  is negative. Further the inequation $8a(a-1) >0$ implies that
  $\sqrt{b ( b+ 8a(a-1))} > b$ and so
  $b (1-b) + (b-1)\sqrt{b ( b+ 8a(a-1))} < 0$, because $(b-1)$ is negative.
  We conclude that $\lambda_1$ is positive. It follows from \cite[Theorem
  2.15 (i)]{dla:qualplan} that the stationary point is a saddle point.
\end{description}

\begin{lemma}\label{lem:s34}
  In the situation of Lemma~\ref{lem:stationarypoints}
  (\ref{lem:stationarypointsCase2}), we have the following dynamical
  behaviour at the stationary points $s_3$ and $s_4$:
  \begin{enumerate}[itemindent=20pt,nosep]
  \item For $b > 1 > a$ and $| 8a(a-1) | > b$ each point is a strong
    attracting focus.
  \item For $b > 1 > a$ and $| 8a(a-1) | \leq b$ each point is an
    attracting node.
  \item For $a > 1 > b$ each point is a saddle point.
  \end{enumerate}
\end{lemma}

\subsection{The Infinite Stationary Points}\label{sec:infstat}

In this section we determine the local phase portrait at the infinite
stationary points using the Poincar\'e compactification.  The basic idea of
the Poincar\'e compactification is a specific projection of a given planar
vector field to the northern and southern hemisphere of the unit sphere and
to extend this vector field to the equator.  The dynamics at the points of
the equator represent then the dynamics at the infinite points of the
planar system. The construction of the Poincar\'e compactification can be
found in \cite[Section 5.1 \& 5.2]{dla:qualplan}.  Throughout this section
we use the notation and the formulas given there.  Our computations are
done in the charts $(U_1,\phi_1)$ and $(U_2,\phi_2)$ representing the front
and right hemisphere and the charts $(V_1,\psi_1)$ and $(V_2,\psi_2)$
covering the back and left hemisphere.

Obviously, the maximal degree of the polynomial right hand sides $P$ and
$Q$ of the system (\ref{eq:dynsys}) is three and we can write them as the
sums of homogeneous polynomials $P=P_2 +P_3$ and $Q=Q_2 +Q_3$ with
$P_2 = xy$, $P_3 = -a(x^3 + x y^2)$, $Q_2 = by^2 + (b-1)x^2 $ and
$Q_3 =-b(y^3+yx^2)$.  Substituting the homogeneous polynomials $P_3$ and
$Q_3$ in the corresponding formula, we obtain that a point $(u,0)$ of
$\mathbb{S}^1 \cap ( U_1 \cup V_1)$ is stationary if and only if
\[
F(u) = -b(u^3+u) - u (-a ( 1+u^2))=(a-b)u (u^2+1)=0.
\]
We conclude that for $a \neq b$ we have exactly one real stationary point
on the $u$-axis, namely $(u,0)=(0,0)$. By contrast, for $a=b$ every point
on the $u$-axis is stationary. Moreover, a point $(u,0)$ of
$\mathbb{S}^1 \cap ( U_2 \cup V_2)$ is stationary if and only if
\[
G(u)=a(u^3+u) - u(-b(1+u^2))= (b-a)u(u^2+1)=0.
\]
As above for $a\neq b$ the only real stationary point on the $u$-axis is
$(u,0)=(0,0)$ and in case $a=b$ every point there is stationary.  We
determine the dynamics near the isolated stationary point $(u,0)=(0,0)$ in
the different charts. Substituting the above homogeneous polynomials as
well as $F'(u)$ and $G'(u)$ in the formulas for the Jacobians we obtain the
matrices
\[
\begin{pmatrix}
    a-b & b-1  \\
    0   & a 
    \end{pmatrix}\qquad \text{resp.}\qquad
\begin{pmatrix}
    b-a & 0 \\
    0 & b
\end{pmatrix}
\]    
which represent the cases if $(0,0)$ belongs to $U_1 \cup V_1$ or
$U_2 \cup V_2$, respectively.  Because of $a\neq b$, the stationary point
$(0,0)$ is in both cases a hyperbolic. In case of the charts $(U_1,\phi_1)$
and $(V_1,\psi_1)$, it is for $b>a$ a saddle point and for $a>b$ a
repelling node by \cite[Theorem 2.15]{dla:qualplan}.  In the charts
$(U_2,\phi_2)$ and $(V_2,\psi_2)$ we obtain by the same theorem that it is
for $b>a$ a repelling node and for $a>b$ a saddle.

As mentioned above, for $a=b$ every point on the $u$-axis in either chart
is stationary and we obtain for the corresponding Jacobians
\begin{displaymath}
\begin{pmatrix}
    0 & (b-1)(u^2+1) \\
    0 & b(u^2+1)
\end{pmatrix}\qquad \text{resp.}\qquad
\begin{pmatrix}
    0 & (1-b)u(u^2+1) \\
    0 & b(u^2+1)
\end{pmatrix}
\end{displaymath}
for points on $U_1\cup V_1$ or $U_2\cup V_2$, respectively.  Obviously, the
centre manifold is given by the complete $u$-axis, as it consists of
stationary points, and is therefore uniquely determined by
\cite[Cor.~3.3]{js:cm}.  Furthermore, in both cases the second eigenvalue
is positive and thus a unique one-dimensional unstable manifold exists
\cite[Thm.~3.2.1]{gh:osc}. This implies that in this case no trajectory can
have an $\omega$-limit point at infinity, but that each point there is the
$\alpha$-limit point of a unique trajectory, namely its unstable manifold.

\begin{lemma}\label{lem:inf}
  \begin{enumerate}[itemindent=20pt,nosep]
  \item If $a=b$, then every point at infinity is a stationary point with
    one outgoing trajectory.
  \item If $a\neq b$, then the points at infinity in positive and negative
    $x$- and $y$-direction are stationary points.
    \begin{enumerate}[leftmargin=15pt,itemindent=25pt,nosep]
    \item If $b>a$, then the stationary points at infinity in positive and
      negative $x$-direction are saddle points and the stationary points at
      infinity in positive and negative $y$-direction are repelling nodes.
    \item If $a>b$, then the stationary points at infinity in positive and
      negative $x$-direction are repelling nodes and the stationary points
      at infinity in positive and negative $y$-direction are saddle points.
    \end{enumerate}
  \end{enumerate}
\end{lemma}

\subsection{Putting Everything Together}

Finally, we collect all our results. As discussed in more detail in the
next section, we show for each arising case a phase portrait in the
appendix.

\begin{theorem}\label{thm:main}
  \begin{enumerate}[itemindent=20pt,nosep]
  \item[(1)] The dynamics of the system (\ref{eq:dynsys}) for the parameter
    values $a=b=1$ is as follows: The finite stationary points are the
    points on the circle $x^2+(y-\tfrac{1}{2})^2=\tfrac{1}{4}$.  At
    infinity, every point is a stationary point. The orbits are straight
    lines starting at infinity in direction to the origin (see Fig.~10(a)).
  \item[(2)] The dynamics of the system (\ref{eq:dynsys}) for the parameter
    values as in (2) of Lemma~\ref{lem:stationarypoints} is:
    \begin{enumerate}[leftmargin=20pt,itemindent=15pt,nosep]
    \item For $b<1<a$ the stationary point $s_1$ has dynamics as in (1) of
      Lemma~\ref{lem:stationaryPoint_s1} and the stationary point $s_2$ is
      an attracting node.  The points $s_3$ and $s_4$ are saddle
      points. The points at infinity in positive and negative $x$-direction
      are repelling nodes and the ones in positive and negative
      $y$-direction are saddle points (see Fig.~9(a)).
    \item For $a<1<b$ and $| 8a(a-1) |>b$ the stationary point $s_1$ has
      dynamics as in (2) of Lemma~\ref{lem:stationaryPoint_s1} and the
      stationary point $s_2$ is a saddle point.  The stationary points
      $s_3$ and $s_4$ are strong attracting foci.  The points at infinity
      in positive and negative $x$-direction are saddle points and the ones
      in positive and negative $y$-direction are repelling nodes (see
      Fig.~11(a)).
    \item For $a<1<b$ and $| 8a(a-1) |\leq b$ the stationary point $s_1$
      has dynamics as in (2) of Lemma~\ref{lem:stationaryPoint_s1} and the
      stationary point $s_2$ is a saddle point.  The stationary points
      $s_3$ and $s_4$ are attracting nodes.  The points at infinity in
      positive and negative $x$-direction are saddle points and the ones in
      positive and negative $y$-direction are repelling nodes (see
      Fig.~11(b)).
    \end{enumerate}
  \item[(3)] The dynamics of the system (\ref{eq:dynsys}) for the parameter
    values as in (3) of Lemma~\ref{lem:stationarypoints} is:
    \begin{enumerate}[leftmargin=20pt,itemindent=15pt,nosep]
    \item For $a=1$ and $1<b$ the stationary point $s_1$ has dynamics as in
      (2) of Lemma~\ref{lem:stationaryPoint_s1} and the stationary point
      $s_2$ is an attracting node. The points at infinity in positive and
      negative $x$-direction are saddle points and the ones in positive and
      negative $y$-direction are repelling nodes (see Fig.~11(c)).
    \item For $a=1$ and $b<1$ the stationary point $s_1$ has dynamics as in
      (1) of Lemma~\ref{lem:stationaryPoint_s1} and the stationary point
      $s_2$ is a saddle point where the $x$-direction is repelling and the
      $y$-direction attracting.  The points at infinity in positive and
      negative $x$-direction are repelling nodes and the ones in positive
      and negative $y$-direction are saddle points (see Fig.~9(b)).
    \item For $1<a$, $1<b$ and $a=b$ the stationary point $s_1$ has
      dynamics as in (2) of Lemma~\ref{lem:stationaryPoint_s1} and the
      stationary point $s_2$ is an attracting node.  At infinity every
      point is a stationary point (see Fig.~11(e)).
    \item For $1<a$, $1<b$ and $a<b$ the stationary point $s_1$ has
      dynamics as in (2) of Lemma~\ref{lem:stationaryPoint_s1} and the
      stationary point $s_2$ is an attracting node. The points at infinity
      in positive and negative $x$-direction are saddle points and the ones
      in positive and negative $y$-direction are repelling nodes (see
      Fig.~11(d)).
    \item For $1<a$, $1<b$ and $a>b$ the stationary point $s_1$ has
      dynamics as in (2) of Lemma~\ref{lem:stationaryPoint_s1} and the
      stationary point $s_2$ is an attracting node. The points at infinity
      in positive and negative $x$-direction are repelling nodes and the
      ones in positive and negative $y$-direction are saddle points (see
      Fig.~11(f)).
    \item For $a<1$, $b<1$ and $a=b$ the stationary point $s_1$ has
      dynamics as in (1) of Lemma~\ref{lem:stationaryPoint_s1} and the
      stationary point $s_2$ is a saddle where the $x$-direction is
      repelling and the $y$-direction is attracting.  At infinity every
      point is a stationary point (see Fig.~9(d)).
    \item For $a<1$, $b<1$ and $a<b$ the stationary point $s_1$ has
      dynamics as in (1) of Lemma~\ref{lem:stationaryPoint_s1} and the
      stationary point $s_2$ is a saddle where the $x$-direction is
      repelling and the $y$-direction is attracting. The points at infinity
      in positive and negative $x$-direction are saddle points and the ones
      in positive and negative $y$-direction are repelling nodes (see
      Fig.~9(e)).
    \item For $a<1$, $b<1$ and $a>b$ the stationary point $s_1$ has
      dynamics as in (1) of Lemma~\ref{lem:stationaryPoint_s1} and the
      stationary point $s_2$ is a saddle where the $x$-direction is
      repelling and the $y$-direction is attracting. The points at infinity
      in positive and negative $x$-direction are repelling nodes and the
      ones in positive and negative $y$-direction are saddle points (see
      Fig.~ 9(c)).
    \item For $b=1$ and $a<1/2$ the stationary point $s_1$ has dynamics as
      in (3a) of Lemma~\ref{lem:stationaryPoint_s1} and the stationary
      point $s_2$ is a saddle where the $x$-direction is repelling and the
      $y$-direction is attracting.  The points at infinity in positive and
      negative $x$-direction are saddle points and the ones in positive and
      negative $y$-direction are repelling nodes (see Fig.~10(b)).
    \item For $b=1$ and $a=1/2$ the stationary point $s_1$ has dynamics as
      in (3b) of Lemma~\ref{lem:stationaryPoint_s1} and the stationary
      point $s_2$ is is a saddle where the $x$-direction is repelling and
      the $y$-direction is attracting.  The points at infinity in positive
      and negative $x$-direction are saddle points and the ones in positive
      and negative $y$-direction are repelling nodes (see Fig.~10(c)).
    \item For $b=1$ and $1/2 < a < 1$ the stationary point $s_1$ has
      dynamics as in (3c) of Lemma~\ref{lem:stationaryPoint_s1} and the
      stationary point $s_2$ is is a saddle where the $x$-direction is
      repelling and the $y$-direction is attracting.  The points at
      infinity in positive and negative $x$-direction are saddle points and
      the ones in positive and negative $y$-direction are repelling nodes
      (see Fig.~10(d)).
    \item For $b=1$ and $1<a$ the stationary point $s_1$ has dynamics as in
      (3d) of Lemma~\ref{lem:stationaryPoint_s1} and the stationary point
      $s_2$ and the stationary point $s_2$ is an attracting node. The
      points at infinity in positive and negative $x$-direction are
      repelling nodes and the ones in positive and negative $y$-direction
      are saddle points (see Fig.~10(e)).
    \end{enumerate}
  \end{enumerate}
\end{theorem}

\section{Phase Portraits and Bifurcations of the CDK System}
\label{sec:ppbif}

\begin{figure}[ht]
  \centering
  \includegraphics[width=0.5\textwidth]{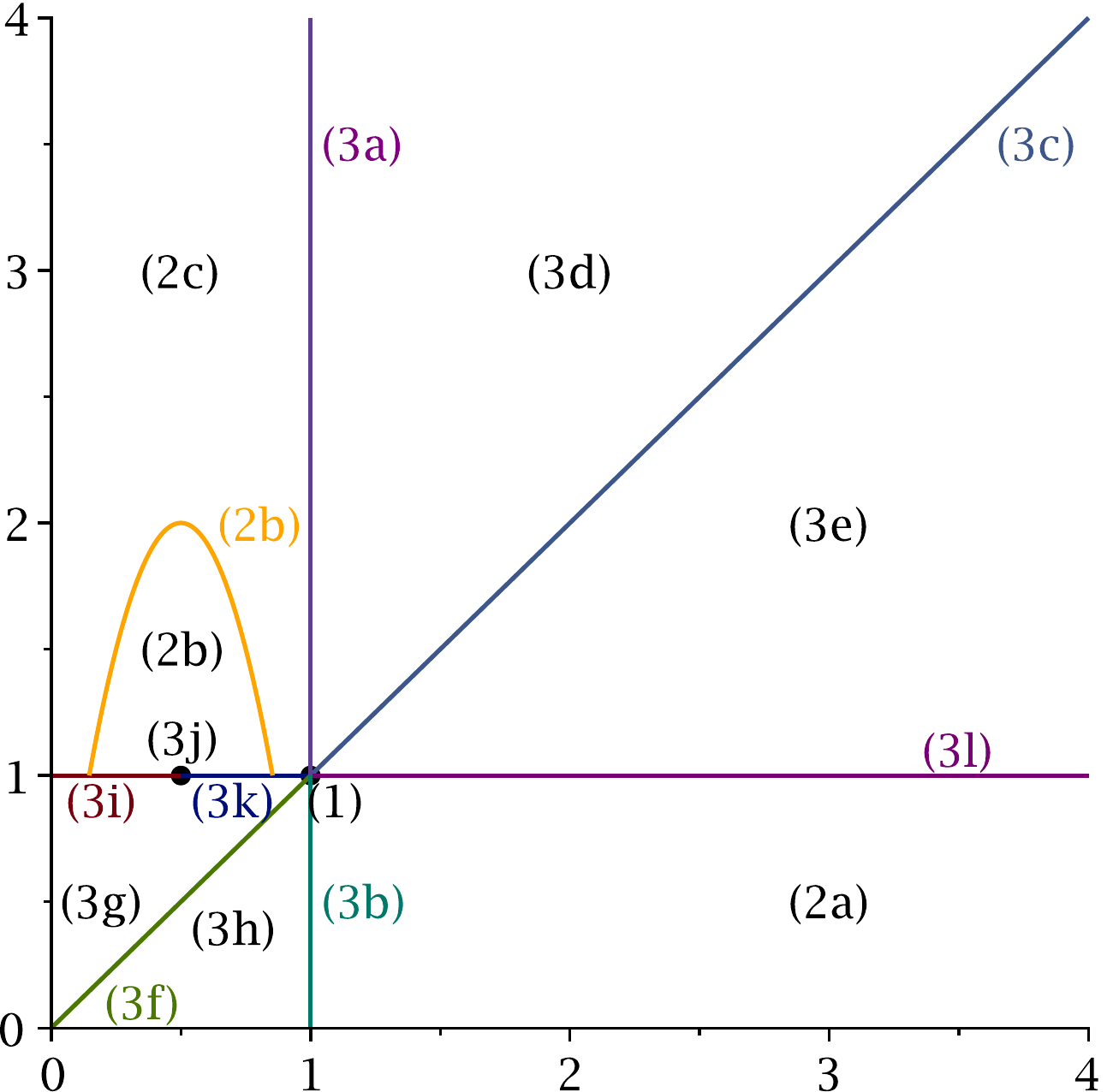}
  \caption{The different cases in the positive orthant of the parameter
    plane}
  \label{fig:para}
\end{figure}

The above analysis of the stationary points distinguishes 16 different
regions in the positive quadrant of the parameter plane.  These are shown
in Fig.~\ref{fig:para} with $a$ on the horizontal axis and $b$ on the
vertical one.  The description in the various regions corresponds to the
numbering of the different cases in Theorem \ref{thm:main}.  One
immediately sees that the main case distinctions arise whenever one of the
parameter crosses the value $1$.

Figs.~\ref{fig:ppbl1}-\ref{fig:ppbg1} in the appendix contain for each
arising case a typical phase portrait.  These phase portraits have been
generated with the help of P4 (see \cite[Chapt.~9]{dla:qualplan} for more
information on this programme).  The phase portraits are shown on the
Poincar\'e disc, i.\,e.\ the infinite plane $\RR^{2}$ is mapped into a
finite disc (see \cite[Chapt.~5]{dla:qualplan} for more details about the
used transformation), as this allows to depict also the stationary points
at infinity.  In the plots, curves in green or cyan consist entirely of
stationary points; curves in red or blue are separatrices.  A green square
signals a saddle, a blue square a stable node and a red square an unstable
node.  Finally, a blue diamond marks a stable strong focus, a blue triangle
a semi-hyperbolic stable node, a green triangle a semi-hyperbolic saddle
and a black x a non-elementary stationary point.

The changes of the phase portraits concern several aspects.  We will mainly
consider the number of stationary points and their types, the number of
attractors and the existence of homoclinic orbits.  The latter point is
controlled exclusively by the parameter $b$.  Therefore we organise our
discussion by the value of $b$.  We will not consider all 16 cases
separately, as sometimes the differences are fairly minor.  In particular,
in most cases the relative size of $a$ and $b$ only affects the stationary
points at infinity.  For $a=b$, all points at infinity are stationary.
Otherwise, there are two stationary points at infinity (antipodes on the
boundary of the Poincar\'e disc can be identified): a saddle and an
unstable node which switch their positions depending on whether $a<b$ or
$b<a$.  We always find either one or two ``almost attractors''.  These are
stationary points and almost all trajectories end in them, the sole
exceptions being trajectories lying in the stable manifolds of saddle
points.

We begin with the case $b<1$ shown in Fig.~\ref{fig:ppbl1}.  A key feature
of the phase portraits is the existence of the two elliptic sectors at the
origin with infinitely many homoclinic orbits.  One effect of them is that
the index of the vector field \eqref{eq:Y} at the origin $s_{1}$ is $2$.
For $a>1$, $s_{2}$ is a stable node and there are the two additional finite
stationary points $s_{3/4}$ which are saddle points (see
Fig.~\ref{fig:dix12}).  We have here two ``almost attractors'', namely
$s_{1}$ and $s_{2}$.  The stable manifolds of the saddle points $s_{3/4}$
separate their basins of attraction.  Furthermore, in this case, infinitely
many heteroclinic orbits connect $s_{1}$ with $s_{2}$.  When $a$ approaches
$1$, then $s_{3/4}$ move towards $s_{2}$.  When $b$ approaches $1$, then
they move towards $s_{1}$.  If $a$ and $b$ simultaneously approach $1$,
then $s_{3}$ and $s_{4}$ move both towards some points on the circle of
stationary points arising for $a=b=1$ (see below) depending on the precise
relationship between $a$ and $b$ in the limit process.  For $a=1$,
$s_{3/4}$ have coalesced with $s_{2}$ to a semihyperbolic saddle and the
heteroclinic orbits connecting $s_{1}$ with $s_{2}$ have all disappeared
except for the one representing a part of the stable manifold of the
saddle.  The separatrix forming the boundary of the elliptic sectors is now
a centre manifold of the semihyperbolic saddle (not shown in
Fig.~\ref{fig:dix13}).  For $a\leq1$, the origin $s_{1}$ is the sole
``almost attractor'' of the system (Figs.~\ref{fig:dix13}-\ref{fig:dix16}).

Now we consider the ``border case'' $b=1$ shown in Fig.~\ref{fig:ppbe1}.
Here we find a very special situation, if also $a=1$ (see
Fig.~\ref{fig:dix01}).  Then the two components of the vector field
\eqref{eq:Y} have a common factor leading to infinitely many finite
stationary points, namely all points on the circle $x^{2}+(y-1/2)^{2}=1/4$
(note that $s_{1}$ and $s_{2}$ both lie on this circle).  As always for
$a=b$, we also find infinitely many stationary points at infinity.  All
trajectories starting at a point $(x_{0},y_{0})$ with $y_{0}\leq0$ or at a
point inside this circle end in $s_{1}$.  All other trajectories connect a
stationary point at infinity with a stationary point on the circle.  If
$a<1$, then we find again the two elliptic sectors at the origin and thus
again the index of the vector field \eqref{eq:Y} at $s_{1}$ is $2$.  The
three cases shown in Fig.~\ref{fig:dix02}-\ref{fig:dix04} differ only in
the computational analysis of the non-elementary stationary point $s_{1}$
which in two cases yields two additional separatrices.  However, each of
these separatrices only separates two parabolic sectors and thus they have
no real influence on the qualitative form of the phase portrait.  We always
have $s_{1}$ as sole ``almost attractor'' of the system.  For $a>1$ (see
Fig.~\ref{fig:dix05}), $s_{2}$ becomes a stable node attracting all
trajectories starting in the upper half plane.  All other trajectories
continue to end in $s_{1}$.  The elliptic sectors become parabolic and
contain infinitely many heteroclinic orbits connecting $s_{1}$ with
$s_{2}$.

Finally, we consider the case $b>1$ shown in Fig.~\ref{fig:ppbg1}.  Here no
elliptic sectors exists and thus no homoclinic orbits.  Furthermore, now
the index of the vector field \eqref{eq:Y} at the origin is $0$.  If $a<1$,
then there are two additional stationary points $s_{3}$ and $s_{4}$ (see
Figs.~\ref{fig:dix06} and \ref{fig:dix07}). Depending on the relative size
of $a$ and $b$, they are either stable focii or stable nodes.  In
particular, they are the only attractors of the system.  For $a$ oder $b$
(or both) approaching $1$, we find the same behaviour as in the case $b<1$
(except that the approach is now for both parameters from the other side of
$1$).  For $a=1$, $s_{3/4}$ have coalesced with $s_{2}$ to a semihyperbolic
stable node (Fig.~\ref{fig:dix08}) and whenever $a\geq1$, then $s_{2}$ is
the sole attractor of the system (Figs.~\ref{fig:dix08}-\ref{fig:dix10}).

It follows from our discussion of the stationary points at infinity in
Section \ref{sec:infstat} that all forward orbits remain bounded, as no
point at infinity has incoming trajectories.  This fact was already
mentioned by \citet{ade:2dchaos}, however, with an erroneous argument.
They claim -- using a Lyapunov type argument -- that any forward orbit
converges into a circle with the radius
$r=\sqrt{\max{\{a^{-2},\tfrac{2-b}{b^2}\}}}$. However, this is not true, as
for large values of $a$ and $b$, this radius can be made arbitrarily small,
nevertheless Lemma~\ref{lem:s2} shows that for these parameter values the
point $(0,1)$ is an attracting node. If one chooses e.g. $a=2.5$, $b=1.9$
(cf.\ Fig.~\ref{fig:dix11}) and considers the point $(0.1,0.9)$ which lies
outside a circle of radius $r=0.4$, then one easily checks that the
derivative $\dot{V}$ of their Lyapunov function is positive at this point
despite the opposite claim of \citet{ade:2dchaos}.

For large values of $b$ and small values of $a$, the radius $r$ determined
by \citet{ade:2dchaos} becomes arbitrarily large. This is indeed necessary,
as for $b>1>a$ one is in case (2) of Lemma \ref{lem:stationarypoints} where
the additional stationary points $s_{3/4}$ exist and these are both
attractive by Lemma \ref{lem:s34}. It is easy to see that for a fixed value
of $b>1$ and for $a\rightarrow0$ the absolute value of the $x$-coordinates
of these points given by \eqref{eq:x12} becomes arbitrarily large.

\section{General Rational Systems} 
\label{sec:ratsys}

So far, we concentrated in this article on the CDK system.  Now we want to
discuss briefly the case of a general planar system with a rational right
hand side.  Let the system be of the form
\begin{equation}\label{eq:ratsys}
  \dot x=\frac{p}{q}\,,\qquad \dot y=\frac{r}{s}
\end{equation}
where $p$, $q$, $r$, $s$ are polynomials in $x$ and $y$ and where we assume
for simplicity that in both equations numerator and denominator are
relatively prime.  Again it is for our purposes more convenient to rewrite
the system in the implicit form
\begin{equation}\label{eq:ratimp}
  q\dot x = p\,,\qquad s\dot y = r\,.
\end{equation}
In the CDK system we are in the special situation that $q=s$ which provides
a minor simplification compared to the general case.

If we introduce $\ell=\lcm{(q,s)}$, a least common multiple of the two
denominators, and write $\ell=q\bar q=s\bar s$, then a straightforward
computation analogous to the one detailed above yields as generator for the
projected Vessiot distribution the vector field
\begin{equation}\label{eq:Yrat}
  Y= \ell\partial_{t} + \bar{q}p\partial_{x} +
                       \bar{s}r\partial_{y}\,.
\end{equation}
This is actually a minimal generator.  Indeed, by definition of a least
common multiple $\gcd{(\bar q,\bar s)}=1$.  Furthermore,
$\gcd{(\ell,\bar{q}p)}=\bar{q}$ and $\gcd{(\ell,\bar{s}r)}=\bar{s}$, since
we assumed that our right hand sides are in reduced form.  Hence
$\gcd{(\ell, \bar{q}p, \bar{s}r)}=1$.

Again one can decouple the $t$-component and thus finds that the
trajectories of (\ref{eq:ratsys}) are the same as the ones of the
polynomial system
\begin{equation}\label{eq:polyrat}
  \dot x = \bar{q}p\,,\qquad
  \dot y = \bar{s}r\,.
\end{equation}
However, because of the decoupling, we can no longer assume that the right
hand sides are relatively prime.  If $d=\gcd{(p,r)}$ is a non-constant
polynomial, then all zeros of $d$ are stationary points leading to a
situation similar to the case $a=b=1$ in the CDK system.  The system
\eqref{eq:polyrat} can now be completely analysed in the same manner, as we
did for the CDK system.

\section{A Logarithmic Variant of the CDK System}

\citet[Sect.~5.3]{jcs:ec} presented further ``chaotic'' two-dimensional
systems.  Some are just minor variations of the CDK system. We have not
analysed them all, but the phase portraits shown in
\cite[Sect.~5.3]{jcs:ec} clearly indicate that most probably all of these
systems also exhibit elliptic sectors instead of chaos.  Finally, Sprott
produced a logarithmic ``chaotic'' system:
\begin{equation}\label{eq:log}
  \dot{x}=\frac{1}{2}\ln{(x^{2})}-y\,,\qquad
  \dot{y}=\frac{1}{2}\ln{(x^{2})}+x\,.
\end{equation}
It is defined only away from the $y$-axis, as for $x=0$ the logarithmic
term becomes infinite. Again it is trivial to show that this system cannot
exhibit any chaotic behaviour.  Multiplication of the right hand side of
\eqref{eq:log} by $x^2$ yields a system which has the same orbits away from
the $y$-axis and which is everywhere in the plane defined and $C^1$. Thus
the Poincar\'e-Bendixson Theorem applies and excludes chaotic
behaviour. For a more detailed analysis of the behaviour near the $y$-axis,
we cannot use the same approach as in the main text. However, it turns out
that very elementary considerations suffice.

\begin{figure}[ht]
  \centering
  \includegraphics[width=0.7\textwidth]{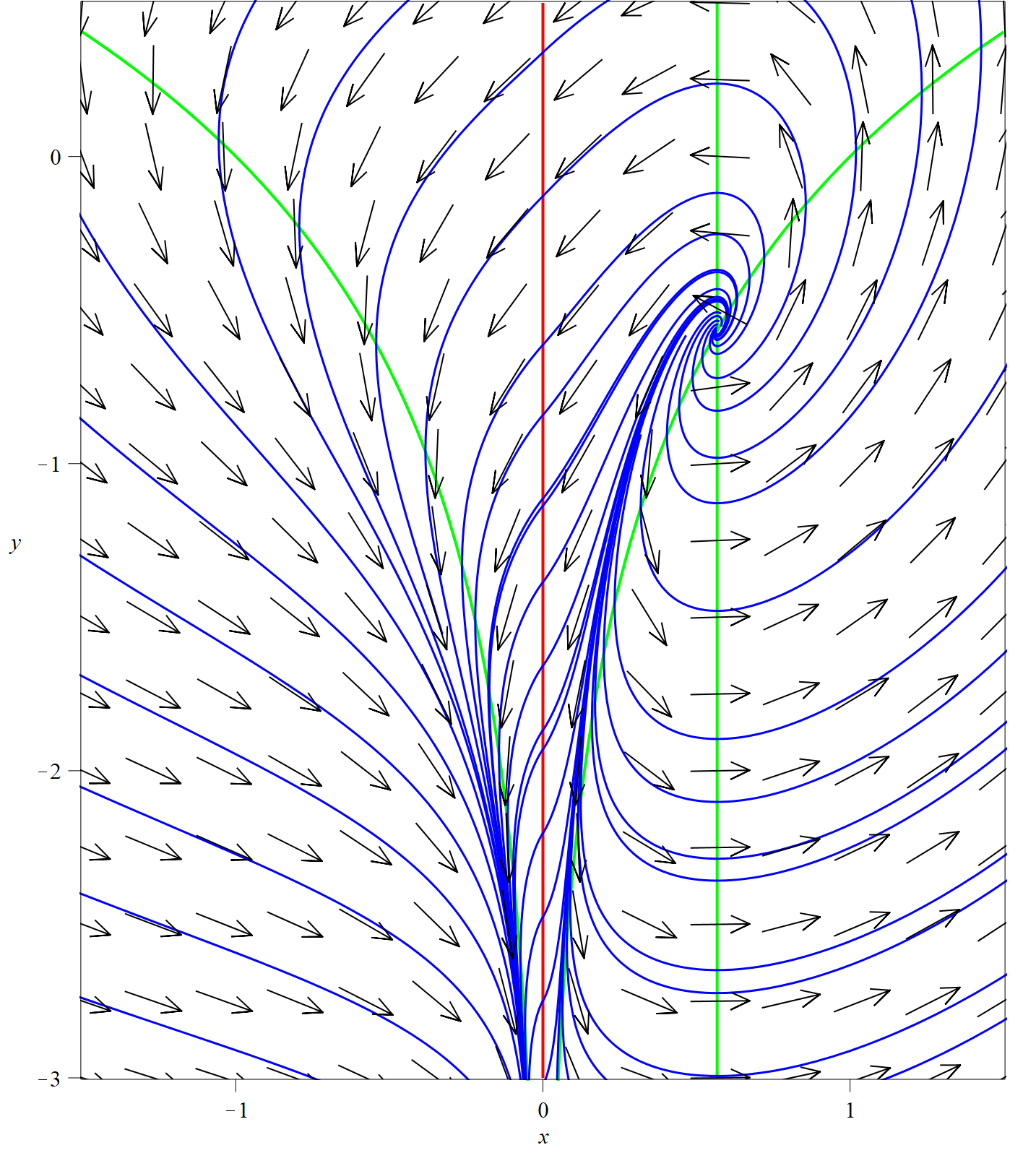}
  \caption{Streamlines and nullclines for system \eqref{eq:log}}
  \label{fig:sprott}
\end{figure}

Figure \ref{fig:sprott} shows some streamlines and the nullclines of system
\eqref{eq:log}.  One clearly sees the one stationary point of the system at
$\bigl(W(1),-W(1)\bigr)\cong(0.56714,-0.56714)$ where $W$ denotes the
Lambert $W$ function. It is an unstable focus. The $y$-nullcline is the
straight line $x=W(1)$, whereas the $x$-nullcline consists of the two
branches of $y=\ln{|x|}$. The streamlines cross the $y$-axis smoothly with
the constant slope $1$, since
\begin{equation}
    \lim_{(x,y)\rightarrow(0,y_0)}{\frac{dy}{dx}}=
    \lim_{(x,y)\rightarrow(0,y_0)}{\frac{1+2x/\ln{x^2}}{1-2y/\ln{x^2}}}=1
\end{equation}
for any value of $y_0$. As the two branches of the $x$-nullcline are
converging towards the negative $y$-axis for $x\rightarrow0$, in the
neighbourhood of a point $(0,y_0)$ with $y_0$ a large negative number the
streamlines are rapidly changing from almost vertical to slope $1$ to again
almost vertical.  The ``chaos'' observed by Sprott is therefore nothing but
a consequence of problems in numerically resolving the dynamics in such a
region with a rapidly changing vector fields.

\section{Conclusions}

In our opinion, the use of the term ``chaos'' in the context of the CDK
system is based on some misconceptions about its mathematical meaning.
Firstly, not every (highly) irregular dynamics is chaotic.  Secondly,
chaotic behaviour is a property of the \emph{exact} trajectories of a
dynamical system and not concerned with numerical difficulties in their
determination or with the effect of noise.  Precise definitions differ from
author to author.  A classical and very popular one goes back to
\citet[Def.~8.5]{rld:chaos}.  In the version given by
\citet[Def.~30.0.2]{sw:dynsyschaos}, one expects the following three
properties from a (compact) chaotic invariant set $\Lambda$ (see also the
discussion in \cite[Def.~7.2.1]{kh:dynsys}): (i) the flow exhibits a
\emph{sensitive dependence on initial conditions} everywhere on $\Lambda$,
(ii) $\Lambda$ is \emph{topologically transitive} and (iii) \emph{periodic
  orbits} are dense in $\Lambda$.  \citet{sw:dynsyschaos} also provides a
series of examples demonstrating the importance of each aspect of these
conditions.  In the case of the CDK system, none of these three properties
is actually given -- even if one restricts to the union of the two elliptic
sectors as a compact invariant set.

In the CDK system for $b<1$ and in Sprott's logarithmic variant one finds
the following situation.  The systems exhibit a sensitive dependence on
initial conditions only in very small regions: for the CDK system, these
regions include in particular that part of the elliptic sectors which is
very close to the origin; in the logarithmic variant, it is a tubular
neighbourhood of the $y$-axis for sufficiently large negative values.  This fact alone would not lead
to the observed very irregular trajectories.  The second key ingredient is
that trajectories leaving the respective sensitive region will reenter it
after a finite time.  In the CDK system, this fact is due to the
ellipticity of the sectors and thus the existence of homoclinic orbits.  In
the logarithmic variant, the unstable focus is mainly responsible for this
behaviour.

If one tries in such a situation to integrate numerically one of the
trajectories starting in the respective sensitive region ``until the end'',
then it is almost unavoidable that at some stage numerical errors lead to a
jump from any incoming trajectory to an outgoing one.  As after a finite
time, this outgoing trajectory will become again an incoming one, the same
will happen repeatedly.  This simple observation explains the highly
irregular behaviour seen in some numerical simulations without any resort
to chaotic dynamics.  One should, however, always keep in mind that the
thus obtained numerical curves do not approximate a \emph{single}
trajectory, but a ``concatenation'' of many different ones (and thus such
curves cannot be used to estimate Lyapunov exponents of single
trajectories).

This phenomenon will appear in any two-dimensional system with a stationary
point possessing at least one elliptic sector (which is possible only at
non-elementary stationary points).  If noise is added to such a system, it
can indeed model highly irregular dynamics without being chaotic and
constructing such a model was the original goal of \citet{cdk:dynmod}.  The
term ``piecewise deterministic dynamics'' coined later by \citet{ddd:pdd}
is for such behaviour much more appropriate than ``chaos''.

The existence of elliptic sectors for the CDK system was already observed
by \citet{ade:2dchaos}.  In contrast to claims by \citet{jcs:ec} and
\citet{ade:2dchaos}, the union of the two elliptic sectors does not define
a two-dimensional attractor, but solely an attracting set.  One of the
simplest definition of an attractor (see e.\,g.\ \cite[Def.~1,
Sect.~3.2]{lp:deds}) says that it is an attracting set which contains a
dense orbit.  Obviously, this is not the case here.  In the relevant region
in the parameter space, the CDK system does not possess at all an attractor
(the origin attracts almost all trajectories in its neighbourhood -- except
the ones starting on the positive $y$-axis), an observation which again
excludes the existence of chaos.

The main contribution of this article is a complete bifurcation analysis of
the CDK system. While \citet{cdk:dynmod} distinguished for their numerical
experiments only four regions in the parameter space, we showed that
including the behaviour at infinity leads actually to sixteen different
regions in the parameter space. If one studies only a finite neighbourhood
of the origin, then essentially one finds only the four cases of
\cite{cdk:dynmod} (plus the boundary cases). The only exception is a subtle
difference in the phase portraits for $b>1$ and $a<1$ where depending on
the value of $a$ the additional stationary points $s_3$ and $s_4$ exhibit
slightly different behaviour. Given the great symmetry of the CDK system
and the fact that physically the parameters $a$, $b$ represent very similar
quantities, it is remarkable that the two parameters possess completely
different roles in the bifurcation analysis. The appearance of the elliptic
sectors is controlled exclusively by $b$, whereas $a$ determines the
behaviour at infinity.

\nonumsection{Acknowledgments} \noindent The authors thank Peter de
Maesschalck (Hasselt University) and Sebastian Walcher (RWTH Aachen) for
many helpful discussions about dynamical systems and Oscar Saleta Reig
(Universitat Aut\`onoma de Barcelona) for helping with P4.  This work was
supported by the bilateral project ANR-17-CE40-0036 and DFG-391322026
SYMBIONT.

\bibliographystyle{ws-ijbc}
\bibliography{Dixon.bib}

\begin{thebibliography}{18}
\newcommand{\enquote}[1]{``#1''}
\providecommand{\natexlab}[1]{#1}
\providecommand{\url}[1]{\texttt{#1}}
\providecommand{\urlprefix}{URL }
\expandafter\ifx\csname urlstyle\endcsname\relax
  \providecommand{\doi}[1]{doi:\discretionary{}{}{}#1}\else
  \providecommand{\doi}{doi:\discretionary{}{}{}\begingroup
  \urlstyle{rm}\Url}\fi

\bibitem[{Alvarez-Ramirez \emph{et~al.}(2005)Alvarez-Ramirez, Delgado-Fernandez
  \& Espinosa-Paredes}]{ade:2dchaos}
Alvarez-Ramirez, J., Delgado-Fernandez, J. \& Espinosa-Paredes, G. [2005]
  \enquote{The origin of a continuous two-dimensional "chaotic" dynamics,}
  \emph{Int.\ J. Bif.\ Chaos} \textbf{15},  3023--3029.

\bibitem[{Arnold(1988)}]{via:geoode}
Arnold, V. [1988] \emph{Geometrical Methods in the Theory of Ordinary
  Differential Equations}, $2^{nd}$ ed., Grundlehren der mathematischen
  Wissenschaften 250 (Springer-Verlag, New York).

\bibitem[{Cummings \emph{et~al.}(1992)Cummings, Dixon \& Kaus}]{cdk:dynmod}
Cummings, F., Dixon, D. \& Kaus, P. [1992] \enquote{Dynamical model of the
  magnetic field of neutron stars,} \emph{Astrophys.\ J.} \textbf{386},
  215--221.

\bibitem[{Devaney(1989)}]{rld:chaos}
Devaney, R. [1989] \emph{An Introduction to Chaotic Dynamical Systems},
  $2^{\text{nd}}$ ed. (Addison-Wesley, Redwood City).

\bibitem[{Dixon(1995)}]{ddd:pdd}
Dixon, D. [1995] \enquote{Piecewise deterministic dynamics from the application
  of noise to singular equations of motion,} \emph{J.\ Phys.\ A: Math.\ Gen.}
  \textbf{28},  5539--5551.

\bibitem[{Dixon \emph{et~al.}(1993)Dixon, Cummings \& Kaus}]{dck:ccd2d}
Dixon, D., Cummings, F. \& Kaus, P. [1993] \enquote{Continuous ``chaotic''
  dynamics in two dimensions,} \emph{Physica D} \textbf{65},  109--116.

\bibitem[{Dumortier \emph{et~al.}(2006)Dumortier, Llibre \&
  Art\'es}]{dla:qualplan}
Dumortier, F., Llibre, J. \& Art\'es, J. [2006] \emph{Qualitative Theory of
  Planar Differential Systems}, Universitext (Springer-Verlag, Berlin).

\bibitem[{Fackerell(1985)}]{edf:vessiot}
Fackerell, E. [1985] \enquote{Isovectors and prolongation structures by
  {V}essiot's vector field formulation of partial differential equations,}
  \emph{Geometric Aspects of the Einstein Equations and Integrable Systems},
  ed. Martini, R., Lecture Notes in Physics~239 (Springer-Verlag, Berlin), pp.
  303--321.

\bibitem[{Guckenheimer \& Holmes(1990)}]{gh:osc}
Guckenheimer, J. \& Holmes, P. [1990] \emph{Nonlinear Oscillations, Dynamical
  Systems, and Bifurcations of Vector Fields}, Applied Mathematical Sciences 42
  (Springer-Verlag, New York).

\bibitem[{Katok \& Hasselblatt(1995)}]{kh:dynsys}
Katok, A. \& Hasselblatt, B. [1995] \emph{Introduction to the Modern Theory of
  Dynamical Systems}, Encyclopedia of Mathematics and its Applications~54
  (Cambridge University Press, Cambridge).

\bibitem[{Perko(2001)}]{lp:deds}
Perko, L. [2001] \emph{Differential Equations and Dynamical Systems}, 3rd ed.,
  Texts in Applied Mathematics~7 (Springer-Verlag).

\bibitem[{Seiler(2010)}]{wms:invol}
Seiler, W. [2010] \emph{Involution --- {T}he Formal Theory of Differential
  Equations and its Applications in Computer Algebra}, Algorithms and
  Computation in Mathematics~24 (Springer-Verlag, Berlin).

\bibitem[{Seiler(2013)}]{wms:singbif}
Seiler, W. [2013] \enquote{Singularities of implicit differential equations and
  static bifurcations,}  \emph{Computer Algebra in Scientific Computing ---
  {CASC} 2013}, eds. Gerdt, V., Koepf, W., Mayr, E. \& Vorozhtsov, E., Lecture
  Notes in Computer Science~8136 (Springer-Verlag, Berlin), pp. 355--368.

\bibitem[{Seiler \& Sei{\ss}(2018)}]{wms:quasilin}
Seiler, W. \& Sei{\ss}, M. [2018] \enquote{Singular initial value problems for
  scalar quasi-linear ordinary differential equations,} Preprint Kassel
  University.

\bibitem[{Sijbrand(1985)}]{js:cm}
Sijbrand, J. [1985] \enquote{Properties of center manifolds,} \emph{Trans.\
  AMS} \textbf{289},  431--469.

\bibitem[{Sprott(2010)}]{jcs:ec}
Sprott, J. [2010] \emph{Elegant Chaos} (World Scientific, Hackensack).

\bibitem[{Vessiot(1924)}]{ves:int}
Vessiot, E. [1924] \enquote{Sur une th\'eorie nouvelle des probl\`emes
  g\'en\'eraux d'int\'e\-gra\-tion,} \emph{Bull.\ Soc.\ Math.\ Fr.}
  \textbf{52},  336--395.

\bibitem[{Wiggins(2003)}]{sw:dynsyschaos}
Wiggins, S. [2003] \emph{Introduction to Applied Nonlinear Dynamical Systems
  and Chaos}, $2^{\text{nd}}$ ed., Text in Applied Mathematics~2
  (Springer-Verlag, New York).

\end{thebibliography}

\nonumsection{Appendix: Phase Portraits for the Various Cases}

For each case arising in Theorem \ref{thm:main}, we have produced a
numerical phase portrait using P4 for one specific set of parameter values
in the corresponding region of the parameter space. A discussion of the
phase portraits and an explanation of the used symbols has already been
given in Section \ref{sec:ppbif}.

\begin{figure}[ht]
  \centering
  \subfigure[\ $a=2.5,\ b=0.5$\label{fig:dix12}]{\includegraphics[width=0.375\textwidth]{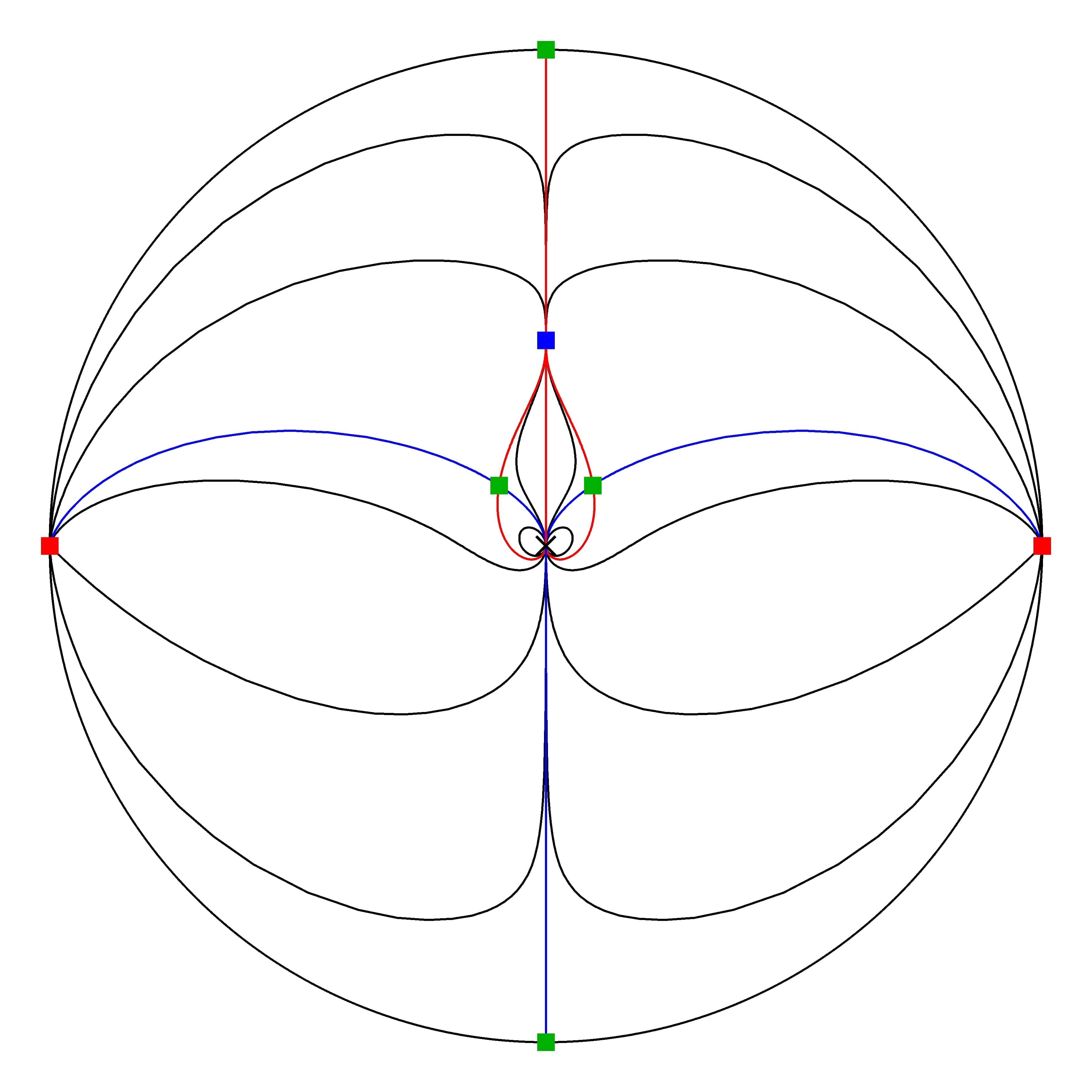}}\quad
  \subfigure[\ $a=1.0,\ b=0.5$\label{fig:dix13}]{\includegraphics[width=0.375\textwidth]{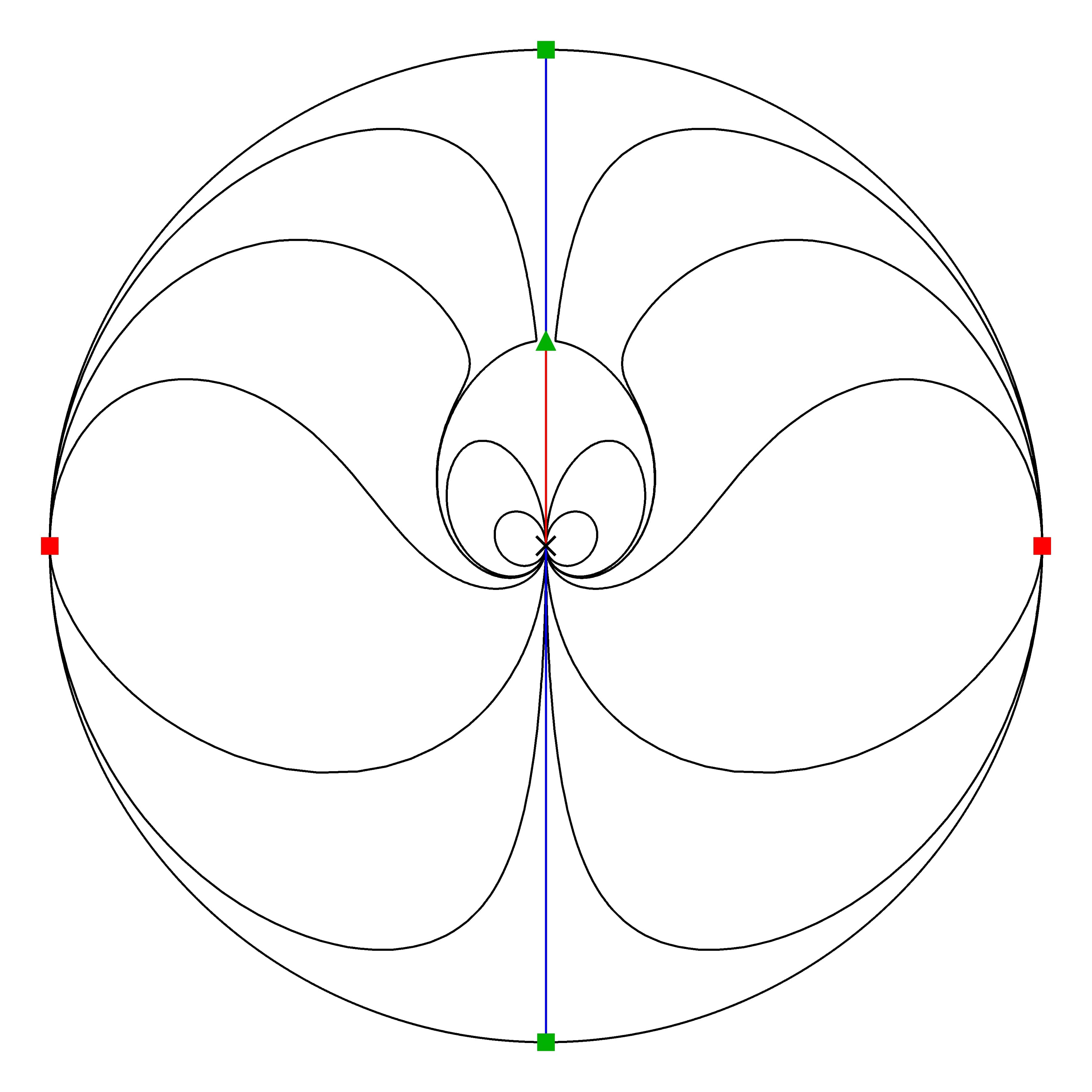}}\\
  \subfigure[\ $a=0.7,\ b=0.5$\label{fig:dix14}]{\includegraphics[width=0.375\textwidth]{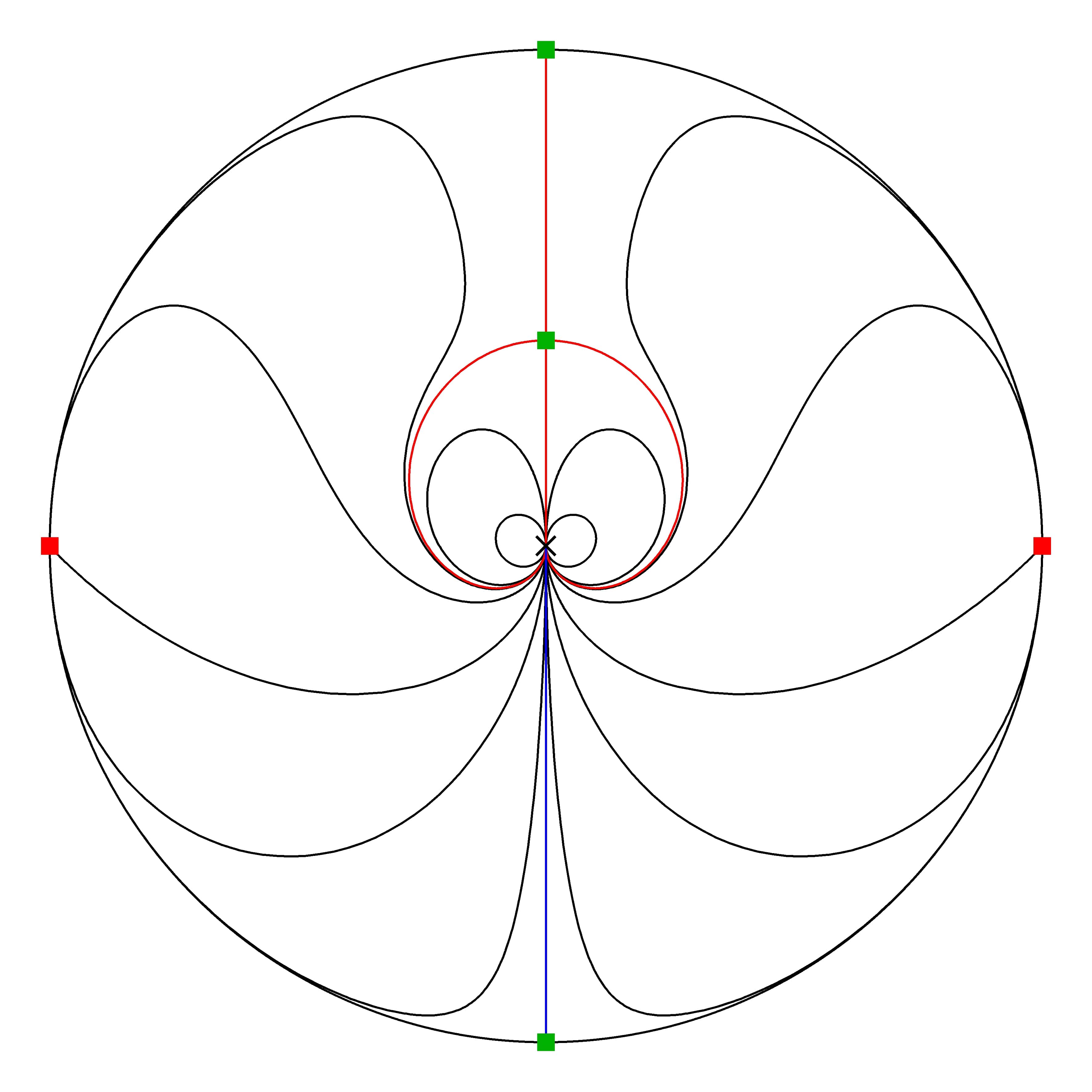}}\quad
  \subfigure[\ $a=0.5,\ b=0.5$\label{fig:dix15}]{\includegraphics[width=0.375\textwidth]{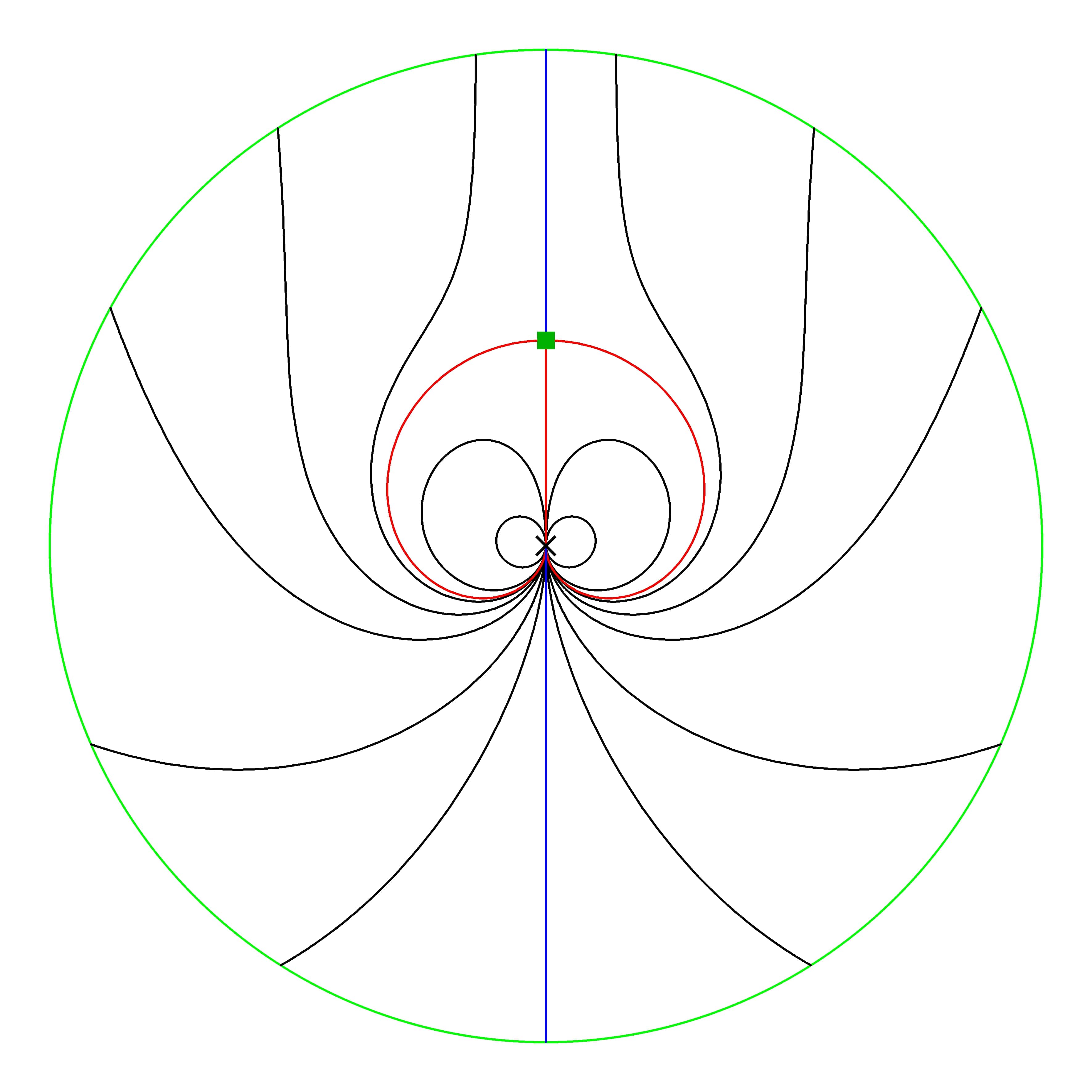}}\\
  \subfigure[\ $a=0.2,\ b=0.5$\label{fig:dix16}]{\includegraphics[width=0.375\textwidth]{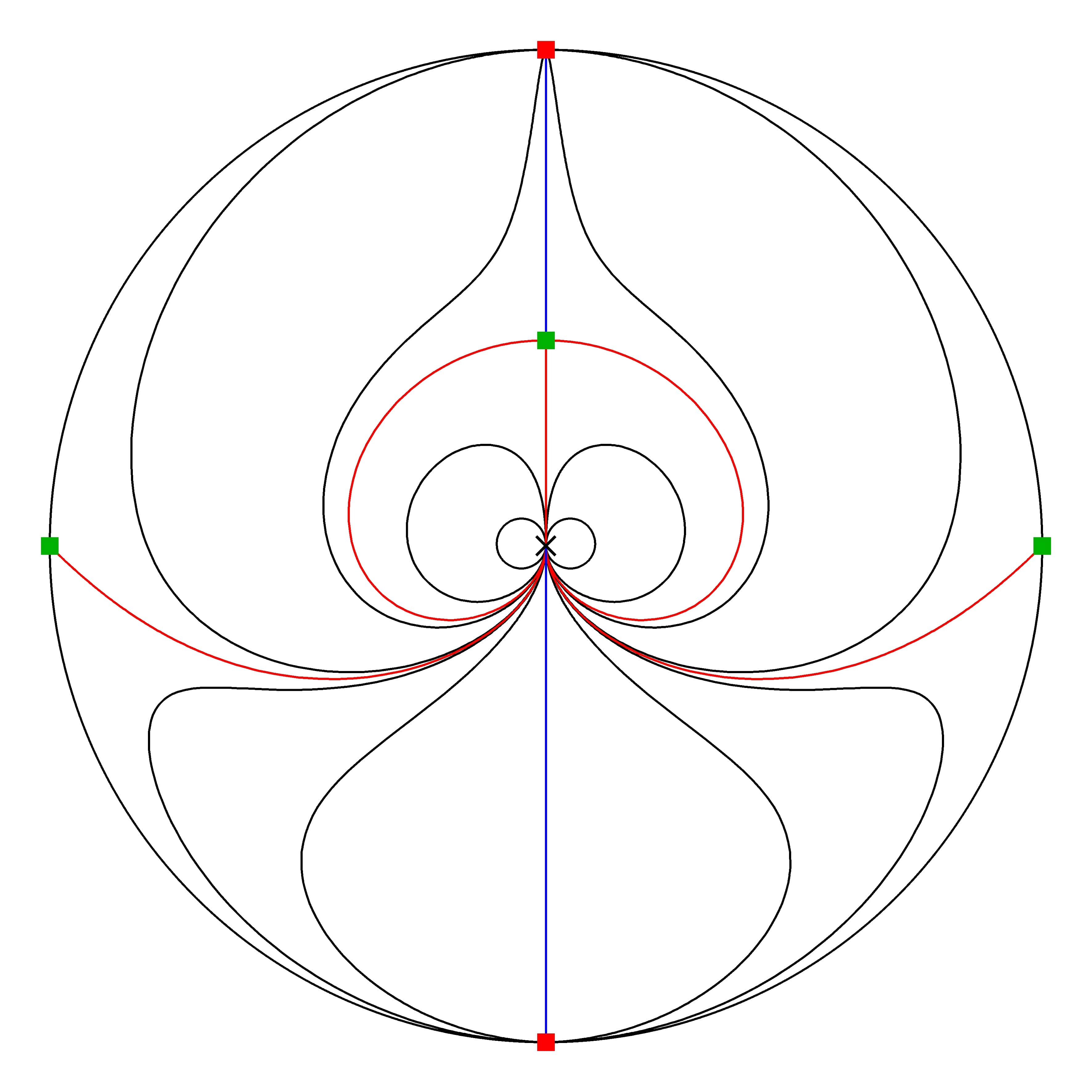}}
  \caption{Phase portraits of the CDK system for $b=0.5$ and different
    values of $a$}\label{fig:ppbl1}
\end{figure}

\begin{figure}[ht]
  \centering
  \subfigure[\ $a=1.0,\ b=1.0$\label{fig:dix01}]{\includegraphics[width=0.375\textwidth]{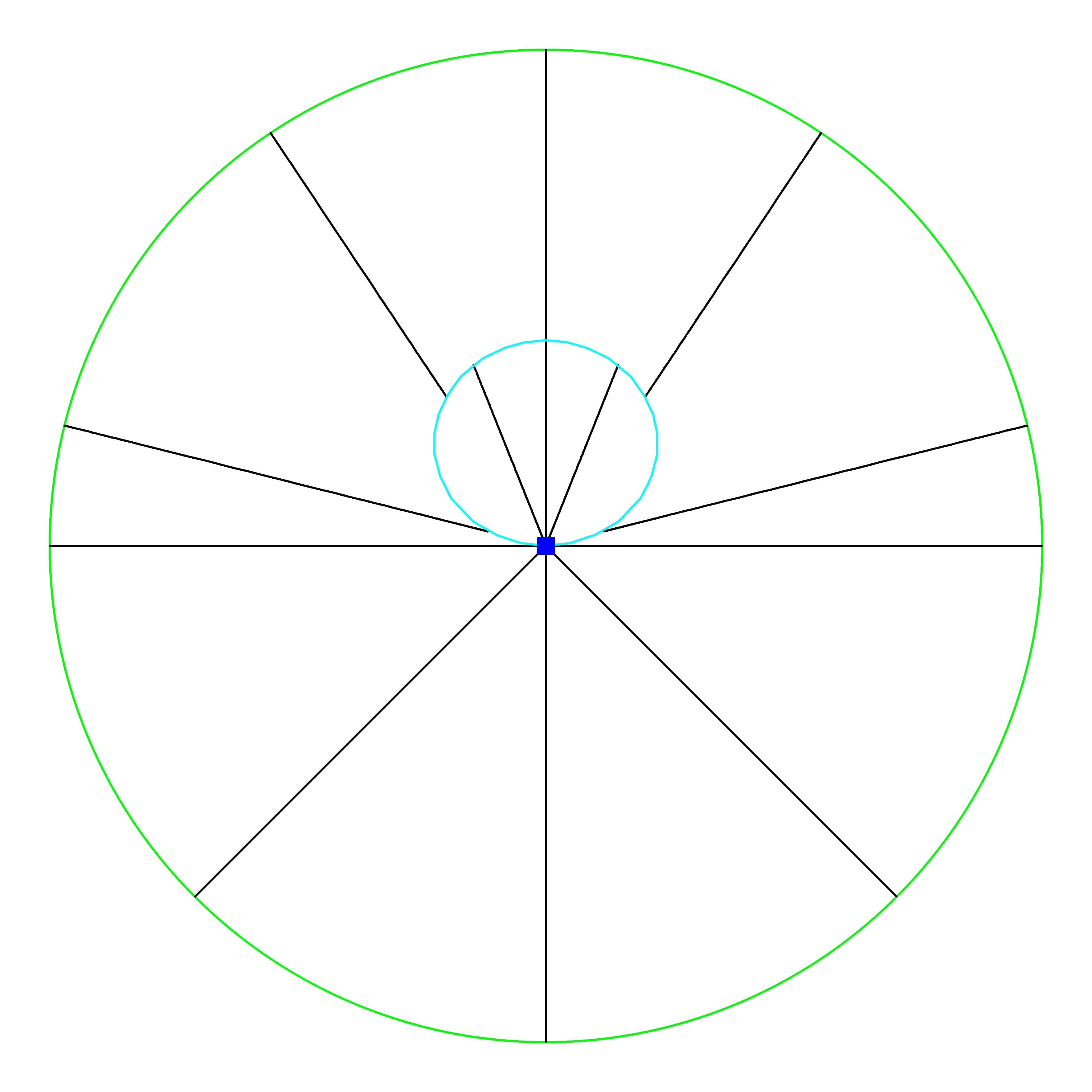}}\quad
  \subfigure[\ $a=0.2,\ b=1.0$\label{fig:dix02}]{\includegraphics[width=0.375\textwidth]{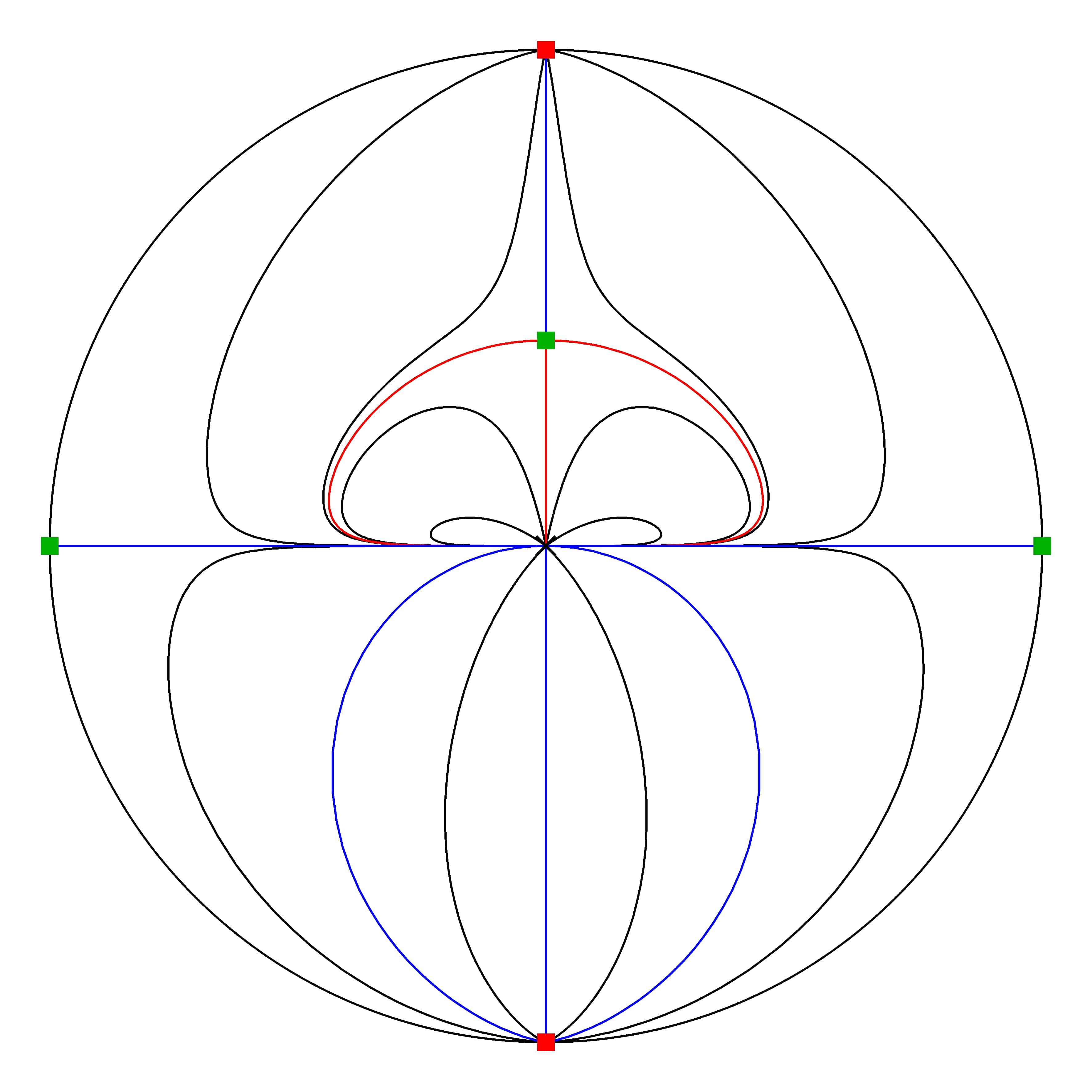}}\\
  \subfigure[\ $a=0.5,\ b=1.0$\label{fig:dix03}]{\includegraphics[width=0.375\textwidth]{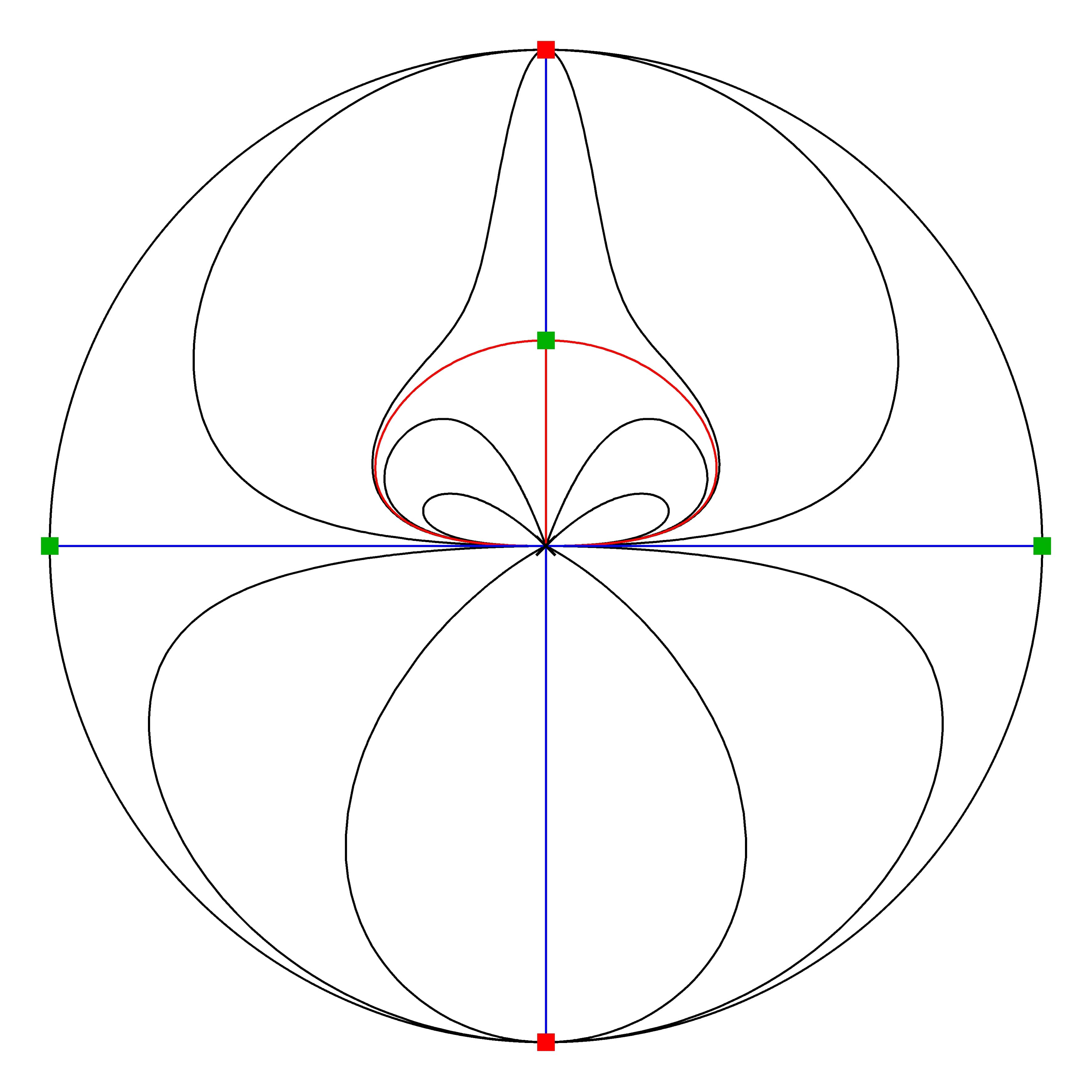}}\quad
  \subfigure[\ $a=0.7,\ b=1.0$\label{fig:dix04}]{\includegraphics[width=0.375\textwidth]{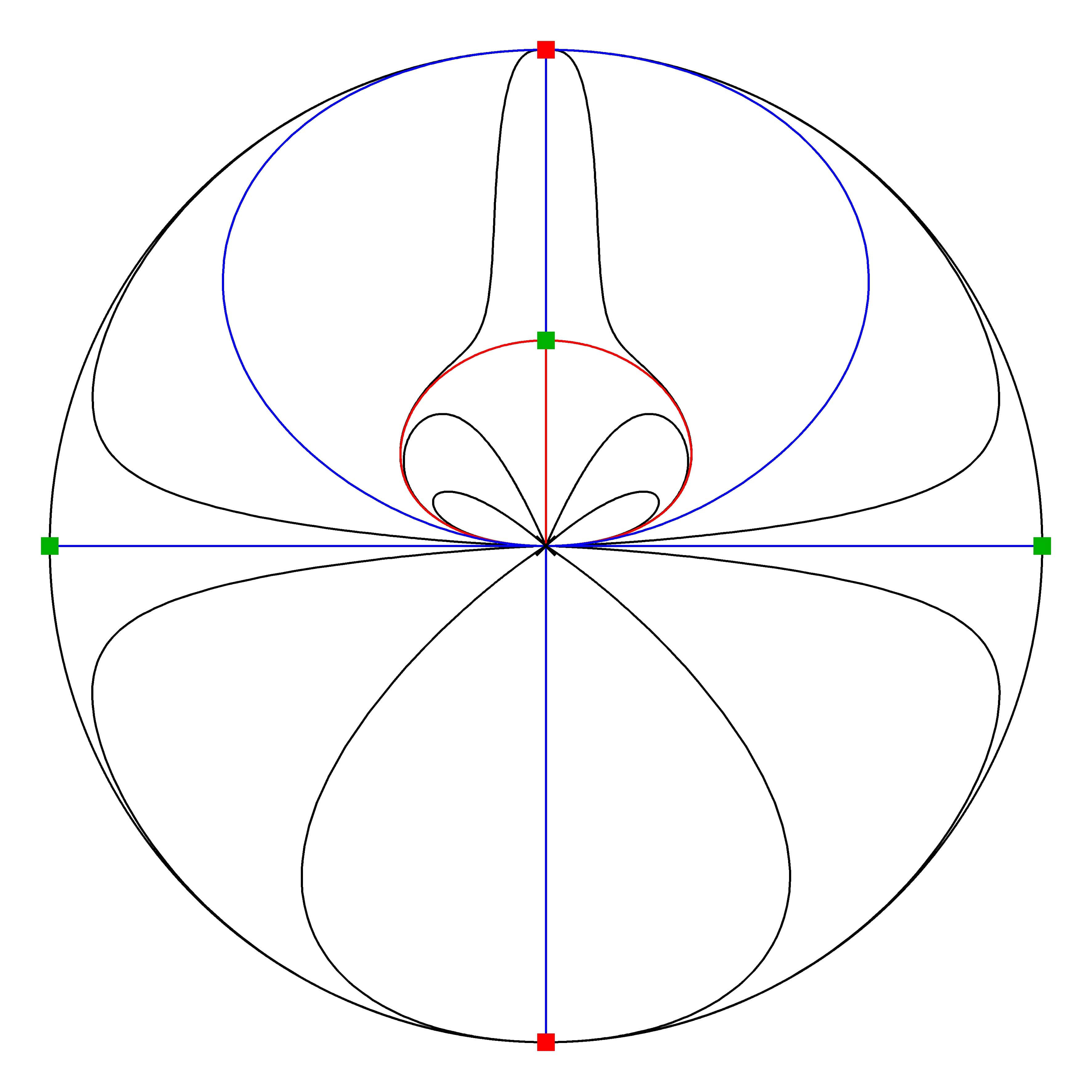}}\\
  \subfigure[\ $a=2.5,\ b=1.0$\label{fig:dix05}]{\includegraphics[width=0.375\textwidth]{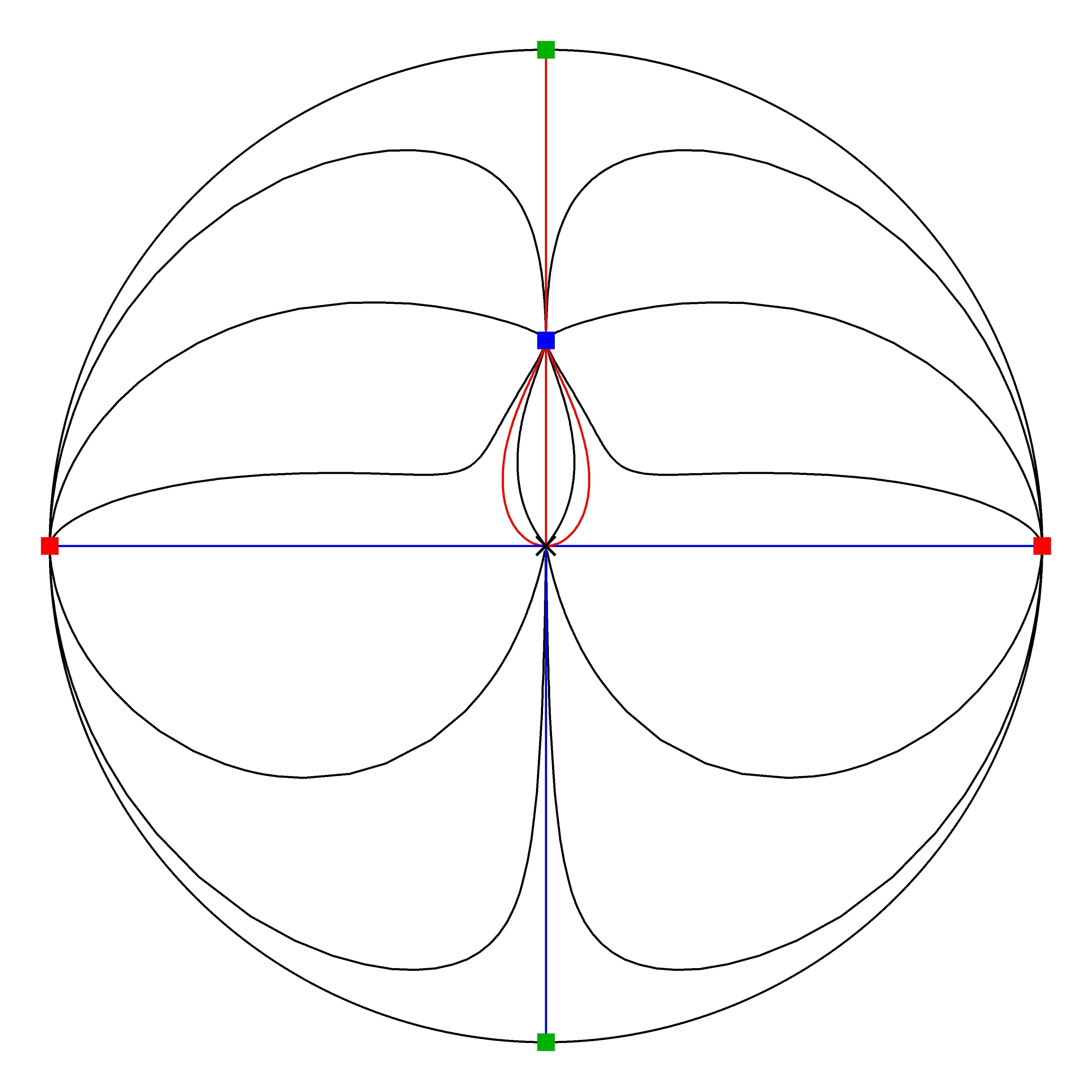}}
  \caption{Phase portraits of the CDK system for $b=1$ and different
    values of $a$}
  \label{fig:ppbe1}
\end{figure}

\begin{figure}[ht]
  \centering
  \subfigure[\ $a=0.5,\ b=1.9$\label{fig:dix06}]{\includegraphics[width=0.375\textwidth]{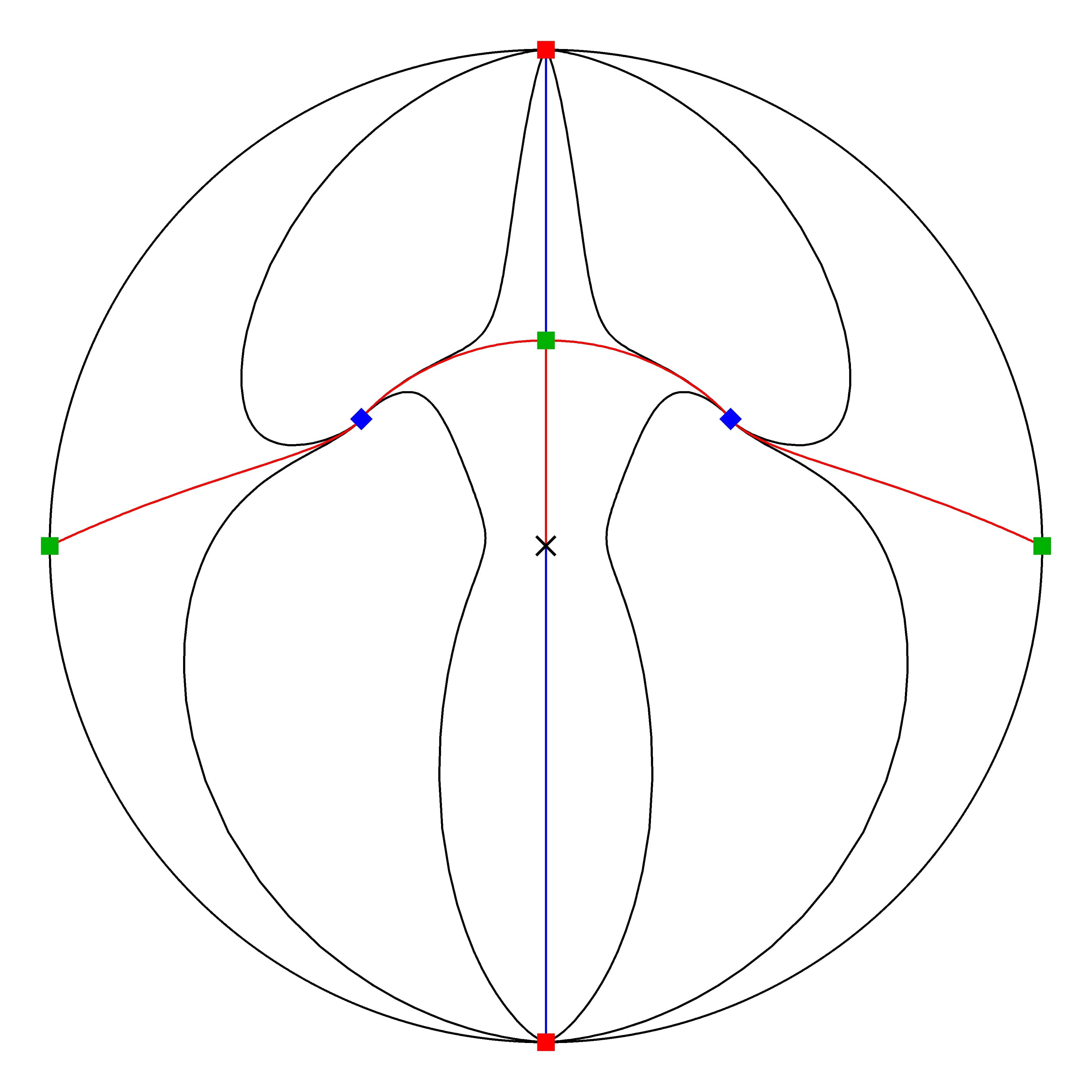}}\quad
  \subfigure[\ $a=0.7,\ b=1.9$\label{fig:dix07}]{\includegraphics[width=0.375\textwidth]{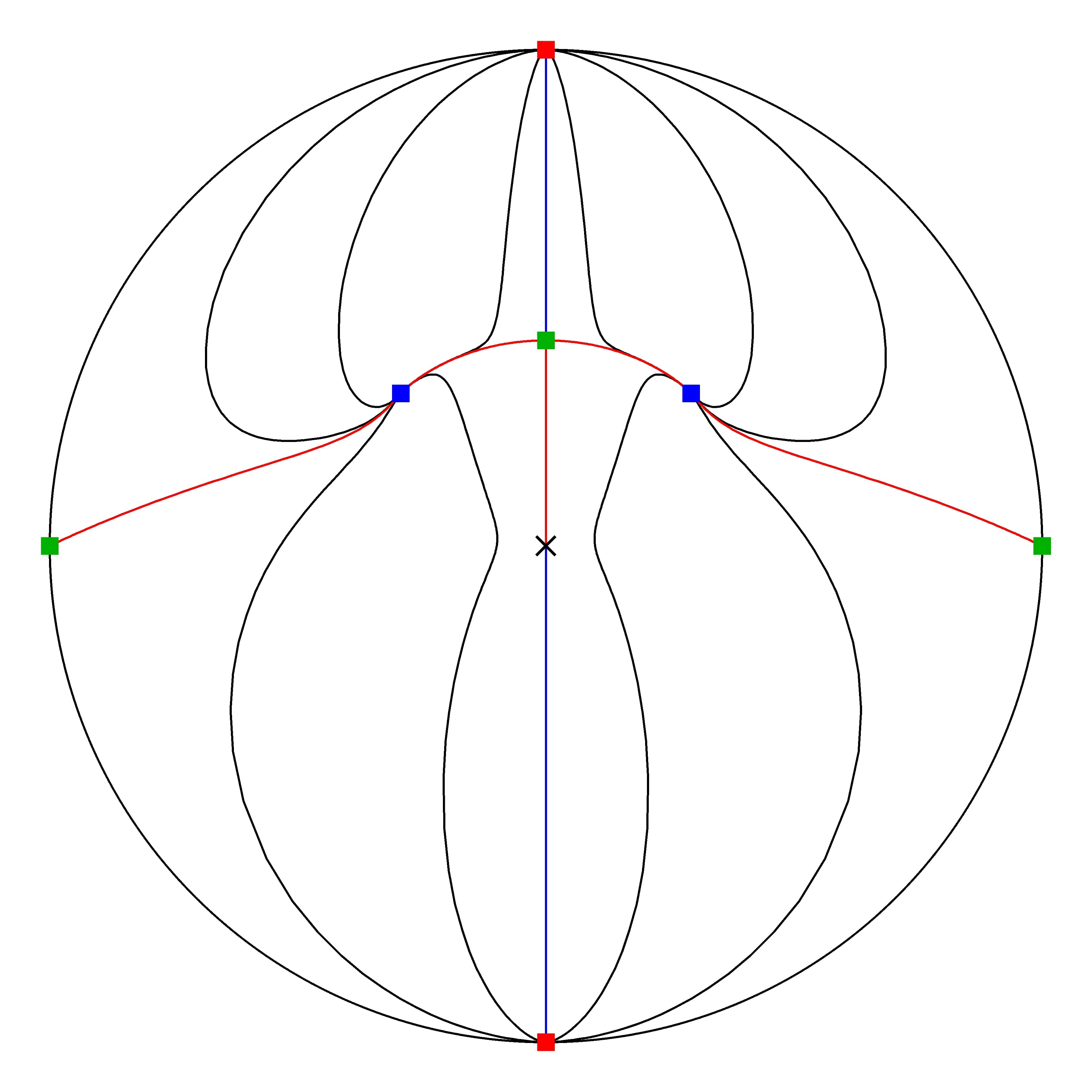}}\\
  \subfigure[\ $a=1.0,\ b=1.9$\label{fig:dix08}]{\includegraphics[width=0.375\textwidth]{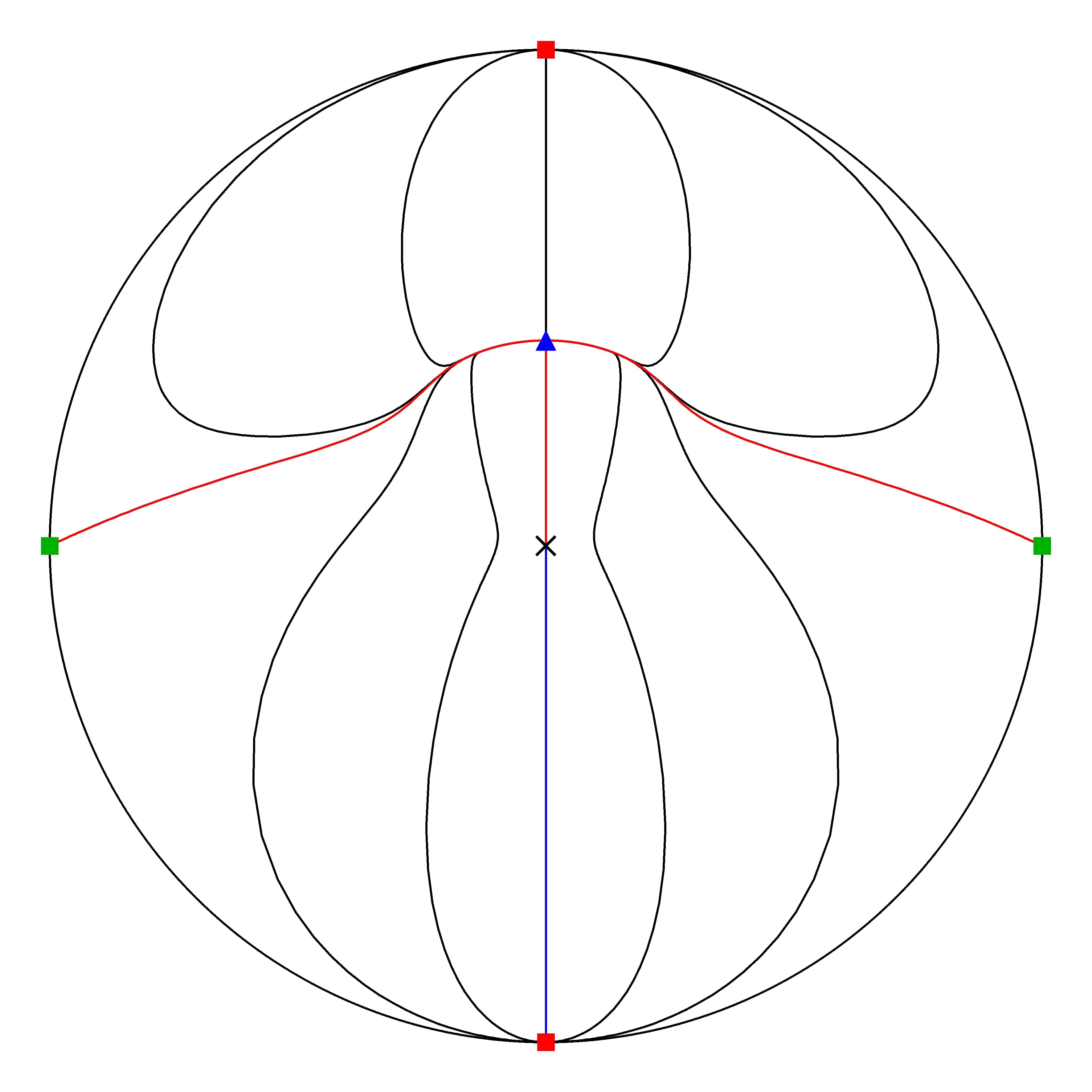}}\quad
  \subfigure[\ $a=1.2,\ b=1.9$\label{fig:dix09}]{\includegraphics[width=0.375\textwidth]{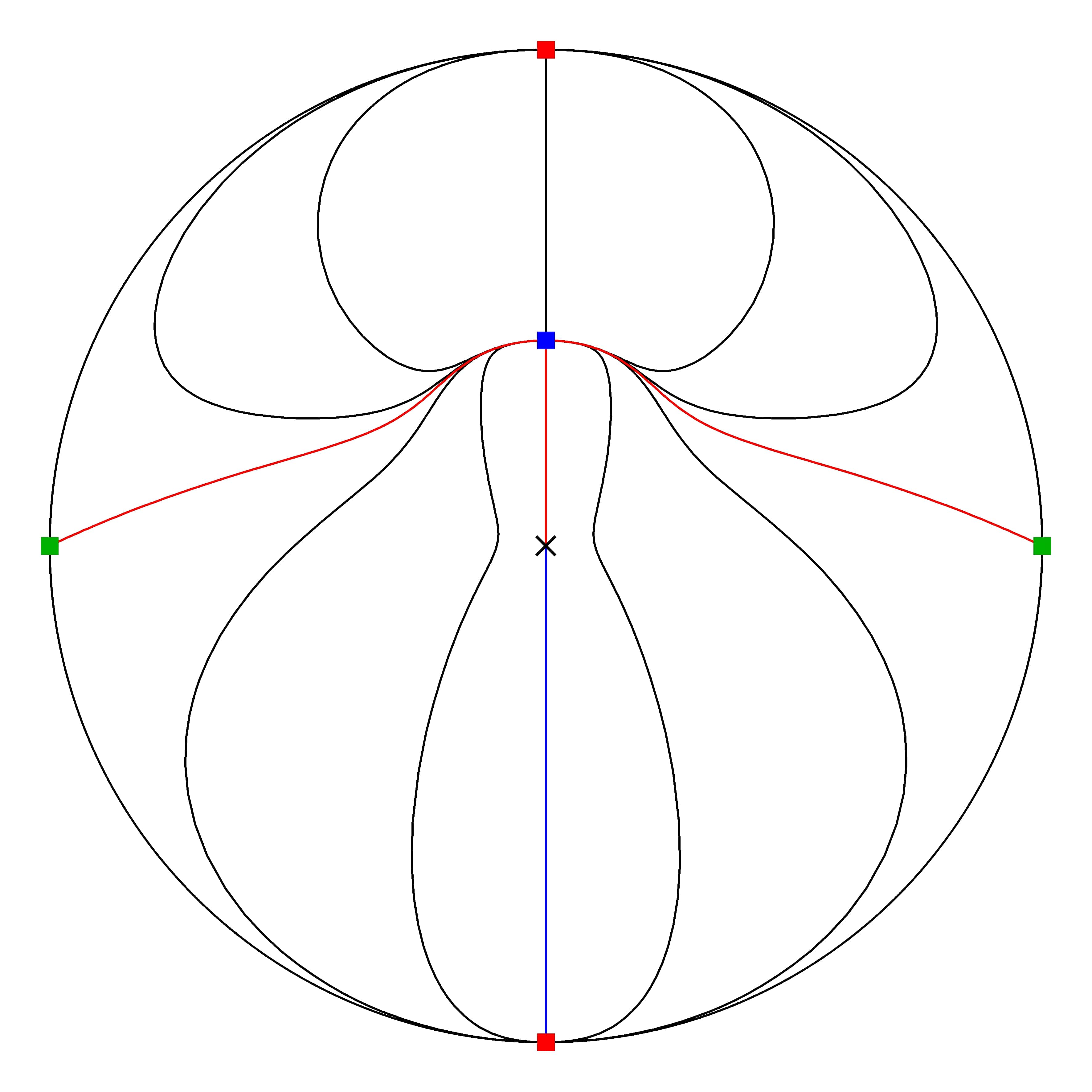}}\\
  \subfigure[\ $a=1.9,\ b=1.9$\label{fig:dix10}]{\includegraphics[width=0.375\textwidth]{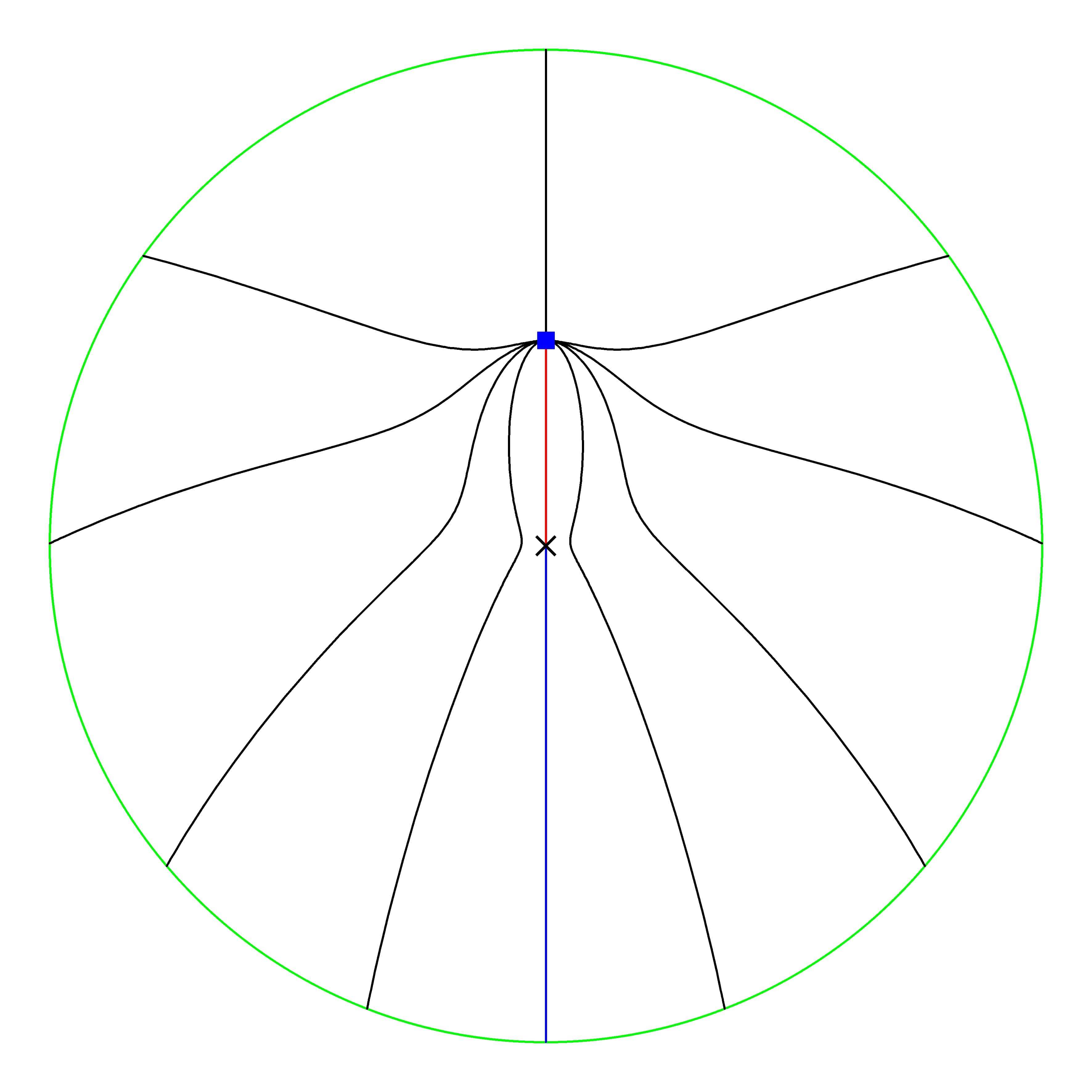}}\quad
  \subfigure[\ $a=2.5,\ b=1.9$\label{fig:dix11}]{\includegraphics[width=0.375\textwidth]{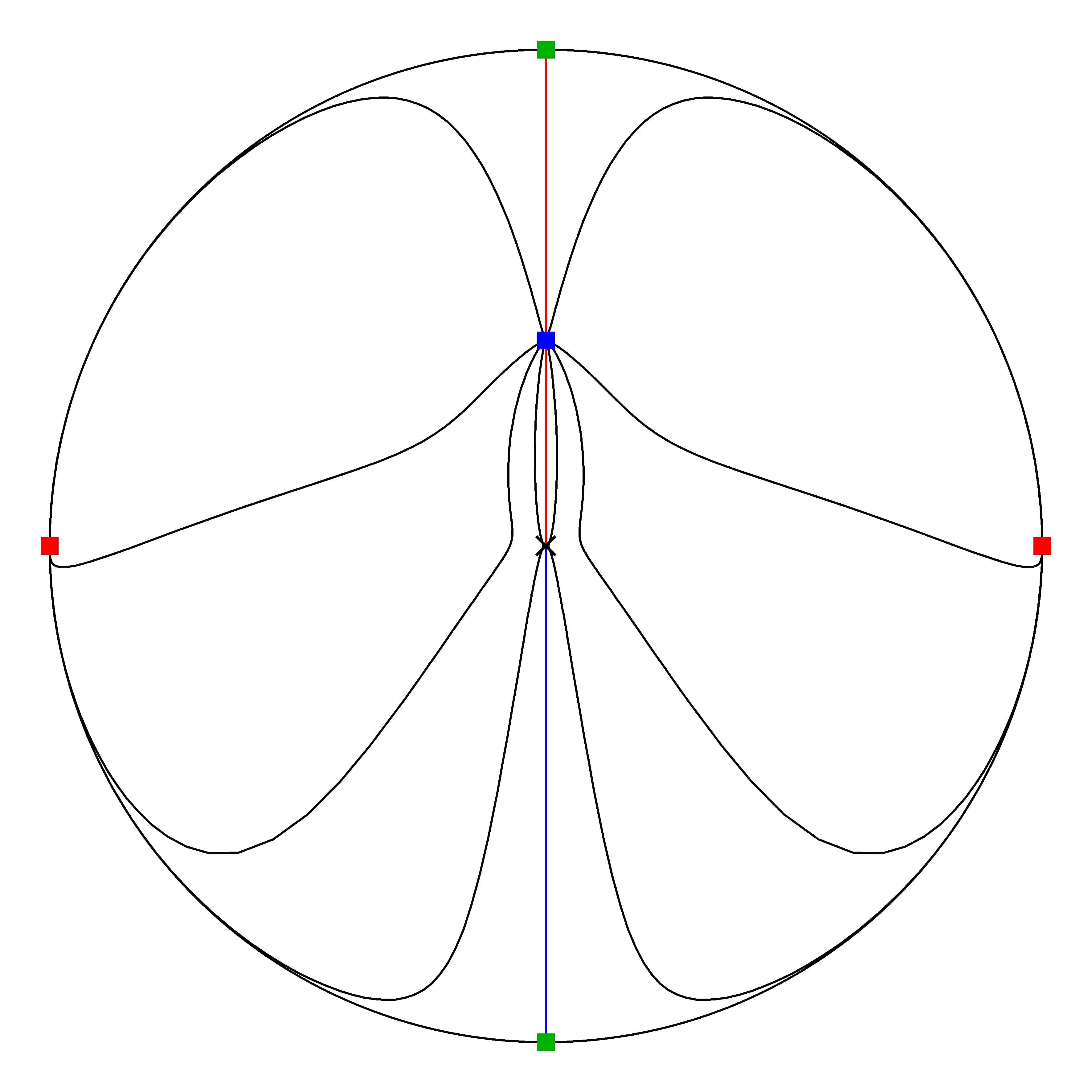}}
  \caption{Phase portraits of the CDK system for $b=1.9$ and different
    values of $a$}\label{fig:ppbg1}
\end{figure}

\end{document}